\documentclass[12pt]{amsart}
\usepackage[a4paper,margin=1in]{geometry}
\usepackage{cmap}
\usepackage[T1]{fontenc}
\usepackage[utf8]{inputenc}
\usepackage{amsmath,amssymb,amsthm,mathtools}
\usepackage{graphicx}
\usepackage{subcaption}
\usepackage{booktabs,array}
\usepackage{enumitem}
\usepackage{xcolor}
\usepackage{tikz-cd}
\usepackage{placeins}
\usepackage[colorlinks=true,linkcolor=red!75!black,citecolor=green!50!black,
            urlcolor=blue!60!black]{hyperref}
\hypersetup{
  pdftitle={Local connectivity of Julia sets of Newton maps of trigonometric
            and hyperbolic functions},
  pdfauthor={Sayantan Chakraborty},
  pdfsubject={MSC 2020: 37F10, 37F20, 30D05, 30D30},
  pdfkeywords={Newton's method, transcendental meromorphic functions,
               Julia set, local connectivity, petals at infinity,
               Whyburn's criterion}}

\graphicspath{{figures/}}
\setlist[enumerate]{label=\textup{(\alph*)},leftmargin=*,itemsep=2pt}

\theoremstyle{plain}
\newtheorem{theorem}{Theorem}[section]
\newtheorem{lemma}[theorem]{Lemma}
\newtheorem{proposition}[theorem]{Proposition}
\newtheorem{corollary}[theorem]{Corollary}
\newtheorem{thmx}{Theorem}

\theoremstyle{definition}
\newtheorem{definition}[theorem]{Definition}
\newtheorem{question}{Question}[section]
\newtheorem{remark}{Remark}[section]

\numberwithin{equation}{section}
\makeatletter
\def\@cite#1#2{{\@citestyle[\citeform{#1}\if@tempswa,\nobreakspace\mbox{#2}\fi]}}
\makeatother

\newcommand{\C}{\mathbb{C}}
\newcommand{\Ch}{\widehat{\mathbb{C}}}
\newcommand{\R}{\mathbb{R}}
\newcommand{\Z}{\mathbb{Z}}
\newcommand{\N}{\mathbb{N}}
\newcommand{\D}{\mathbb{D}}
\newcommand{\B}{\mathcal{B}}
\newcommand{\PS}{\mathcal{P}}
\newcommand{\Pcrit}{\mathcal{P}_{\mathrm{crit}}}
\newcommand{\Pasym}{\mathcal{P}_{\mathrm{asym}}}
\newcommand{\Nsin}{N_{\sin}}
\newcommand{\Ncos}{N_{\cos}}
\newcommand{\Nsinh}{N_{\sinh}}
\newcommand{\Ncosh}{N_{\cosh}}
\DeclareMathOperator{\Arg}{Arg}
\DeclareMathOperator{\Sing}{Sing}
\DeclareMathOperator{\Acc}{Acc}
\DeclareMathOperator{\diam}{diam}
\renewcommand{\Re}{\operatorname{Re}}
\renewcommand{\Im}{\operatorname{Im}}

\begin{document}

\title[Local connectivity of Julia sets of Newton maps]{Local connectivity of Julia sets of Newton maps of trigonometric and hyperbolic functions}

\author{Molla Basir Ahamed$^*$}
\address{Molla Basir Ahamed, Department of Mathematics, Jadavpur University, Kolkata 700032,	\mbox{India}}
\email{mbahamed.math@jadavpuruniversity.in}

\author{Sayantan Chakraborty}
\address{Sayantan Chakraborty, Department of Mathematics, Jadavpur University, Kolkata 700032, \mbox{India}}
\email{chakrabortysayantan51@gmail.com}


\begin{abstract}
	We study the Newton maps of $\sin z$, $\cos z$, $\sinh z$, and $\cosh z$, each of which is transcendental meromorphic and has infinitely many unbounded immediate attracting basins. We first verify directly that, for the cosine Newton map $N_{\cos}(z) = z + \cot z$, the half-planes $\{\pm \operatorname{Im} z > M\}$ are repelling petals at infinity for every $M \ge 2$. These explicit petals allow us to apply the tame-basin criterion of Bara\'nski, Fagella, Jarque, and Karpi\'nska and obtain local connectivity of the spherical boundary of each immediate basin of $N_{\cos}$. We also show that the Newton map of $Ae^{\lambda z} + Be^{-\lambda z}$, where $AB\lambda \neq 0$, is affinely conjugate to the sine Newton map $N_{\sin}(z) = z - \tan z$. The published local-connectivity theorem for $J(N_{\sin})$ therefore applies to this entire family, including the cosine and the two hyperbolic Newton maps. We describe the invariant lines and basin geometry of the four maps and illustrate the corresponding dynamical planes numerically.
\end{abstract}

\maketitle

\section{Introduction}\label{sec:intro}

\subsection{Background and literature review}\label{subsec:background}

Let $f\colon\C\to\Ch$ be a transcendental meromorphic function. We
define its \emph{Fatou set} by
\begin{equation}\label{eq:fatou-set}
	\mathcal{F}(f)=
	\left\{
	z\in\C:
	\begin{array}{l}
		\text{there exists an open neighborhood $U\subset\C$ of $z$ on which}\\
		\text{all iterates $f^n$, $n\geq1$, are defined and the family}\\
		\text{$\{f^n|_U:n\geq1\}$ is normal in the spherical metric}
	\end{array}
	\right\}.
\end{equation}
We call the connected components of $\mathcal{F}(f)$ the \emph{Fatou components}
and define the \emph{Julia set} by
\begin{equation}\label{eq:julia-set}
	\mathcal{J}(f)=\Ch\setminus \mathcal{F}(f).
\end{equation}
In particular, $\infty\in \mathcal{J}(f)$.\vspace{1.2mm}

The iteration theory of rational maps goes back to Fatou and Julia
\cite{Fat19,Jul18}, and Fatou subsequently studied transcendental entire
functions \cite{Fat26}. We refer to \cite{Ber93a,Dom98} for meromorphic
dynamics and to \cite{CG93,Mil06} for rational dynamics. We now recall
the definition of local connectivity.
\begin{definition}[\cite{Mun00}]\label{def:lc}
A topological space $X$ is \emph{locally connected} if for every $x\in X$ and every open neighborhood $V$ of $x$ there is a connected open set $W$ with $x\in W\subset V$.
\end{definition}
Local connectivity of the boundary of a Fatou component is a key property for understanding the dynamics of $f$ on its closure. Indeed, if $U$ is a simply connected invariant Fatou component whose boundary in $\Ch$ has more than one point, then a Riemann map $\varphi\colon\D\to U$ conjugates $f|_U$ to a holomorphic self-map of the unit disc $\D$, and by Carathéodory's theorem (see, \textit{e.g.,} \cite{Mil06}) $\varphi$ extends continuously to $\overline{\D}$ if and only if $\partial U$ is locally connected. Local connectivity of the whole Julia set is stronger: by Torhorst's theorem (Lemma~\ref{lem:torhorst}), if $\mathcal{J}(f)$ is locally connected then the boundary of every Fatou component is locally connected. Conversely, if a Fatou component $U$ is completely invariant then $\partial U=\mathcal{J}(f)$, so the two properties coincide; this happens, for instance, for the basin of infinity of a polynomial with connected Julia set.\vspace{1.2mm}

Whether boundaries of Fatou components, or whole Julia sets, are locally
connected depends strongly on the \emph{singular values} of $f$, \textit{i.e.,} the singularities of the inverse function. For rational maps these are exactly the critical values; for transcendental maps one must also take into account the asymptotic values and the accumulation points of critical and asymptotic values. Two main strategies have been developed. The first uses Yoccoz puzzles, built from external rays and equipotentials, and is mainly adapted to polynomials and certain rational maps \cite{Hub93,Mil00}. The second consists in constructing a conformal metric near $\partial U$ (or near $\mathcal{J}(f)$) with
respect to which $f$ is expanding along sufficiently long orbit segments. For hyperbolic maps one may use the hyperbolic metric of a suitable domain in $\Ch\setminus\overline{\PS(f)}$; for subhyperbolic maps, where preperiodic critical points may lie in the Julia set, an orbifold metric with a discrete singular set is used; and for geometrically finite maps the orbifold metric is combined with suitable metrics near the parabolic points \cite{CG93,DH84,Mil06}. Once such a metric is available, one pulls back a curve close to $\partial U$ by inverse branches of $f$; the resulting sequence of curves converges uniformly, so $\partial U$ is a continuous image of a circle and hence locally connected.\vspace{1.2mm}

Using these methods, local connectivity has been established in many settings. For polynomials, the connected Julia sets of large classes of maps are known to be locally connected: for instance, hyperbolic and subhyperbolic polynomials \cite{DH84,DH85} and, under suitable hypotheses, quadratic polynomials treated by Yoccoz's puzzle method \cite{Hub93,Mil00,Yoc99}, real polynomials
\cite{LvS98}, unicritical polynomials \cite{KL09} and non-renormalizable
polynomials \cite{KvS09}, as well as typical quadratic polynomials with a Siegel disc \cite{PZ04}; we refer to these works for the precise statements, and see also \cite{BH88,McM94,Sul83}. Moreover, every bounded Fatou component of a polynomial other than a Siegel disc is bounded by a Jordan curve \cite{RY08,RY22}. For rational maps, connected Julia sets of (sub)hyperbolic maps \cite{CG93,DH84,Mil06} and, more generally, of geometrically finite maps \cite{TY96} are locally connected, and the same holds for wide classes of Newton maps of polynomials \cite{DS22,Roe08,WYZ23}.\vspace{1.2mm}

For transcendental entire functions the situation changes drastically, since the essential singularity at $\infty$ is also an omitted value. In fact, $\mathcal{J}(f)$ cannot be locally connected if $f$ has an unbounded Fatou component \cite{BD00,Osb13}, and an unbounded invariant Fatou component with locally connected boundary must be a univalent Baker domain \cite{BW91,Kis97}. Informally, the preimages of curves that tend to $\infty$ may accumulate in comb-like structures. A well-known example is $z\mapsto\lambda e^z$, $0<\lambda<1/e$, which has an invariant attracting basin of infinite degree whose boundary has a dense set of accesses to $\infty$ \cite{DG87}. A finite-degree example is $z\mapsto ze^{z+1}$, whose unbounded superattracting basin of the fixed point $-1$ has degree $2$ and a similar behaviour, and whose
boundary is conjectured to contain indecomposable continua of escaping points (see \cite{BFJK25} and compare \cite{FJ23}). On the positive side, a locally connected boundary of a bounded simply connected Fatou component of an entire map is a Jordan curve, by the maximum principle; for hyperbolic entire maps the boundaries of bounded Fatou components are Jordan curves, and even quasicircles \cite{BFR15,Mor99}; and the
Julia set is locally connected for hyperbolic entire maps with only bounded Fatou components and no asymptotic values, provided the number of critical points (counted with multiplicity) in each Fatou component is uniformly bounded \cite{BFR15,BM02,Mor99}. Recall that a transcendental map is \emph{geometrically finite} if (with the notation of Section~\ref{subsec:known}) $\Sing(f)\cap \mathcal{F}(f)$ is compact and $\overline{\PS(f)}\cap \mathcal{J}(f)$ is finite, and \emph{strongly geometrically finite} if, in addition, $\mathcal{J}(f)$ contains no asymptotic values and the local degree of $f$ at points of $\mathcal{J}(f)$ is uniformly bounded \cite{ARS22,MB10}. Alhamed \emph{et al.} \cite{ARS22} proved that for strongly geometrically finite entire maps the boundaries of all bounded periodic Fatou components are Jordan curves and, under additional assumptions analogous to the ones above, that all Fatou components are Jordan domains and the Julia set is locally connected. Expanding metrics for transcendental maps were constructed in \cite{ARS22,BFR15,MB10,PS22}. All of these results concern bounded Fatou components and maps in the Eremenko--Lyubich class $\B$ (see Section~\ref{subsec:known}).\vspace{1.2mm}

For transcendental meromorphic maps the essential singularity need not be an omitted value, and unbounded Fatou components with locally connected boundary do occur. For instance, for $z\mapsto\lambda\tan z$ with $\lambda>1$ the Fatou set consists of the upper and the lower half-plane, which are completely invariant attracting basins, so that $\mathcal{J}(f)=\R\cup\{\infty\}$ is a circle in $\Ch$ \cite{DK89}. Numerical experiments suggest that many meromorphic maps have unbounded Fatou components with locally connected boundaries, even when these
components have the essential singularity on their boundary. The first general results in this direction are due to Barański \emph{et al.} \cite{BFJK25}. Roughly speaking, they show that the boundary of a simply connected invariant attracting basin $U$ is locally connected if $f$ is geometrically finite on a compact part of $U$ and the unbounded part of $U$ is contained in finitely many regions in which $f$ is univalent and behaves like a map near a parabolic fixed point (\emph{repelling petals at infinity}). Their results apply to maps outside class $\B$, possibly with asymptotic values and infinitely many post-singular orbits; we recall them in Section~\ref{subsec:known}.

\subsection{Newton's method.} Let $g$ be a non-constant entire function. Its Newton map
\[
N_g(z)=z-\frac{g(z)}{g'(z)}
\]
is meromorphic, and every simple zero of $g$ is a superattracting fixed point of $N_g$. If $g$ is a polynomial, then $N_g$ is rational: Shishikura proved that its Julia set is connected \cite{Shi09}, Hubbard, Schleicher and Sutherland found finite sets of initial points that detect all roots \cite{HSS01}, and Tan Lei studied cubic Newton maps in detail \cite{Tan97}. If $g$ is transcendental, $N_g$ is in general a transcendental meromorphic function. Bergweiler identified a class of meromorphic functions without wandering domains which contains the Newton maps of certain entire functions \cite{Ber93b}; Mayer and Schleicher proved that the immediate basin of every root is simply connected and unbounded \cite{MS06}; Rückert and Schleicher
studied the components of the complement of such a basin \cite{RS07}; and Barański \emph{et al.} proved that all Fatou components of
$N_g$ are simply connected, so that $J(N_g)$ is connected
\cite{BFJK14,BFJK18}. See also \cite{BFJK17,BFJK20} for accesses to infinity from Fatou components and for the relation between Fatou components and singular values. The simplest transcendental example, $\Nsin(z)=z-\tan z$, is the main example of \cite{BFJK25}; an elementary study of its dynamics, including its symmetries and its behaviour on the real line, is given in \cite{CDBSW24} (see Remark~\ref{rem:realtrace}).

\subsection{Tame attracting basins: what is known}\label{subsec:known}
We use the following notation, as in \cite{BFJK25}. Let $f\colon\C\to\Ch$ be meromorphic. A point $v\in\C$ is an \emph{asymptotic value} of $f$ if there is a curve $\gamma\colon[0,\infty)\to\C$ with $\gamma(t)\to\infty$ and
$f(\gamma(t))\to v$ as $t\to\infty$. We denote by $\Sing(f)$ the set of finite singularities of $f^{-1}$, \textit{i.e.,} the critical and asymptotic values of $f$ together with their accumulation points in $\C$. The \emph{post-singular}, \emph{post-critical} and \emph{post-asymptotic} sets are
\begin{gather*}
\PS(f)=\bigcup_{n=0}^{\infty}f^n(\Sing(f))\; \mbox{and}\;
\Pcrit(f)=\{f^n(v): v \text{ a critical value of } f,\ n\ge0\},\\
\Pasym(f)=\{f^n(v): v \text{ an asymptotic value of } f,\ n\ge0\},
\end{gather*}
where the terms for which $f^n$ is not defined are omitted. For $A\subset\C$, the closure $\overline{A}$ and the set $\Acc(A)$ of accumulation points of $A$ are taken in $\C$. The \emph{Eremenko--Lyubich class} $\B$ consists of the
transcendental meromorphic maps $f$ for which $\Sing(f)$ is bounded, \emph{i.e.,} 
\begin{align*}
	\B = \{ f : \mbox{Sing}( f )\; \mbox{is bounded}\}.
\end{align*} We write $\deg_w F$ for the local degree of a holomorphic map $F$ at $w$, and $\diam_{\mathrm{sph}}$ for the diameter with respect to the spherical metric $2|dz|/(1+|z|^2)$. If $d\rho=\rho|dz|$ is a conformal metric on an open set $V\subset\C$, \textit{i.e.,} $\rho$ is a positive continuous function on $V$, then the derivative of a holomorphic map $F$ with respect to $\rho$ is $|F'(z)|_\rho=|F'(z)|\,\rho(F(z))/\rho(z)$.
\begin{definition}[{\cite[Definition~4.1]{BFJK25}}]\label{def:petal}
Let $P\subset\C$ be an unbounded simply connected domain and let
$G\colon P\to\C$ be a holomorphic map which extends continuously to
$\overline{P}$. We say that $P$ is an \emph{attracting petal of $G$ at
infinity} if
\begin{enumerate}
\item $\overline{G(P)}\subset P$;
\item $\bigcap_{n=0}^{\infty}G^n(P)=\emptyset$;
\item for every compact set $K\subset P$ there exist constants $c_1,c_2>0$, $\delta\in(0,\pi/2)$, $n_0\in\N$ and a non-increasing logarithmically convex function $g\colon(t_0,+\infty)\to\R_+$, $t_0>0$, such that $|z|\in(t_0,+\infty)$ and
\[
c_1\,g(|z|)<|G(z)-z|<c_2\,g(|z|) \; \mbox{and}\;
|\Arg(G(z)-z)-\Arg(z)|<\delta
\]
for every $z\in\bigcup_{n\ge n_0}G^n(K)$.
\end{enumerate}
An unbounded simply connected domain $P\subset\C$ is a \emph{repelling petal at infinity} of a univalent map $F\colon P\to\C$ if $F(P)$ is an attracting petal at infinity of $G=F^{-1}$.
\end{definition}
Here, a function $g$ is logarithmically convex if $\ln g$ is convex; the
constant function $g\equiv1$ is the simplest example. We understand the
difference of arguments in (c) modulo $2\pi$, \textit{i.e.,} as the angle between the vectors $G(z)-z$ and $z$.

\begin{definition}[{\cite[Definition~1.1]{BFJK25}}]\label{def:tame}
An invariant attracting basin $U$ of a transcendental meromorphic map $f$ is \emph{tame at infinity} if there exists a disc $D\subset\C$ such that $\overline{U}\setminus D$ is contained in the union of finitely many repelling petals $P_i$ at infinity for $f$, and $U\cap f(P_i)\cap P_{i'}=\emptyset$ for $i\ne i'$.
\end{definition}

\begin{thmx}[{\cite[Theorem~A]{BFJK25}}]\label{thm:A}
Let $U$ be a simply connected invariant attracting basin of a meromorphic map $f\colon\C\to\Ch$. Assume that the following conditions are satisfied.
\begin{enumerate}
\item The set $\bigl(\overline{\Pasym(f)}\cup\Acc(\Pcrit(f))\bigr)\cap
\overline{U}$ is contained in the union of a compact subset of $U$ and a finite set of parabolic periodic points of $f$ in $\partial U$.
\item There exists a compact set $L\subset U$ such that for every
$z\in\Pcrit(f)\cap(\overline{U}\setminus L)$,
\[
\sup\{\deg_w f^n : w\in f^{-n}(z),\ n>0\}<\infty .
\]
\item $U$ is tame at infinity \textup{(}Definition~\ref{def:tame}\textup{)}.
\end{enumerate}
Then the boundary of $U$ in $\Ch$ is locally connected.
\end{thmx}
Theorem~\ref{thm:A} is a consequence of the following result on the existence of an expanding metric near $\partial U$.
\begin{thmx}[{\cite[Theorem~B]{BFJK25}}]\label{thm:B}
Let $U$ be an invariant simply connected attracting basin of a meromorphic map $f\colon\C\to\Ch$ satisfying the assumptions of Theorem~\ref{thm:A}. Then there exists a simply connected domain $A$ with $\overline{A}\subset U$ such that for every compact set $K\subset U$ one can find a conformal metric $d\varsigma=\varsigma|dz|$ on $U\setminus\overline{A}$ and numbers $b_n\in(0,1)$, $n\in\N$, such that $\sum_{n=1}^{\infty}b_n<\infty$ and
\[
|(f^n)'(z)|_\varsigma>\frac{1}{b_n}
\]
for every $z\in U\setminus\overline{A}$ and $n\in\N$ with $f(z),\dots,f^{n-1}(z)\in U\setminus\overline{A}$ and $f^n(z)\in K\setminus\overline{A}$.
\end{thmx}
Let us briefly describe the construction behind Theorem~\ref{thm:B}; see
\cite[Sections~4--7]{BFJK25} for the details. One considers a domain $V$ close to $\partial U$ and the hyperbolic orbifold $(V,\nu)$, where $\nu(z)$ is the least common multiple of the local degrees $\deg_w f^n$, $w\in f^{-n}(z)$, $n\in\N$; roughly speaking, hypotheses (a) and (b) of Theorem~\ref{thm:A} ensure that $(V,\nu)$ is a hyperbolic orbifold. The orbifold metric, \textit{i.e.,} the push-forward of the hyperbolic metric of $\D$ under the universal branched covering of $(V,\nu)$, is then uniformly expanding on compact subsets (see \cite[Theorem~5.1]{BFJK25} for a general version). Near the parabolic points in $\partial U$ and in the repelling petals at infinity it is replaced by metrics of the form $|dz|/|z-p|^{\alpha}$, $0\le\alpha<1$, and $|dz|/|z|^{\alpha}$, $\alpha>1$, respectively, which are expanding in repelling petals. Gluing these metrics gives $\varsigma$. Finally, Theorem~\ref{thm:A} follows by pulling back a curve in $U$ along inverse branches of $f$ and showing, with the help of $\varsigma$, that the resulting curves converge uniformly in the spherical metric.\vspace{1.2mm}

The main example of \cite{BFJK25} is the following.

\begin{figure}[htbp]
\centering
\begin{subfigure}{0.6\textwidth}
\centering
\includegraphics[width=\textwidth]{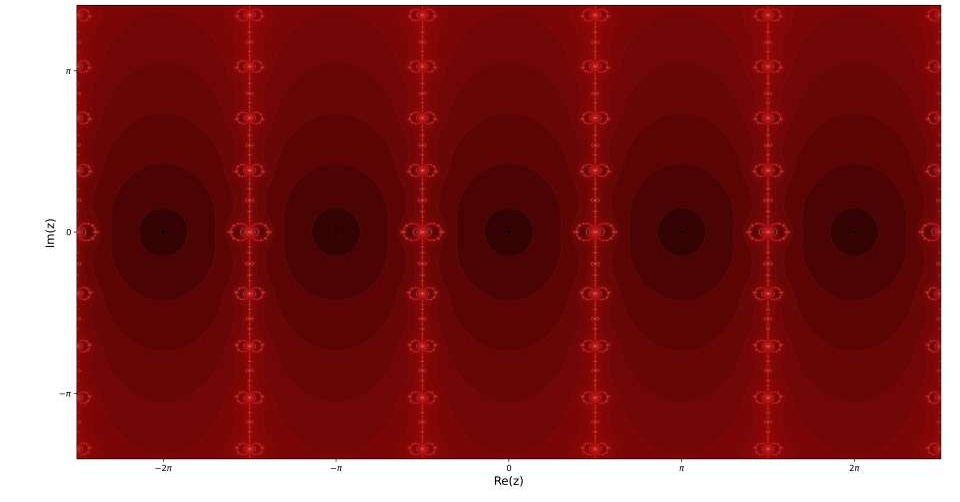}
\caption{The dynamical plane of $\Nsin$.}
\end{subfigure}

\medskip
\begin{subfigure}{0.6\textwidth}
\centering
\includegraphics[width=\textwidth]{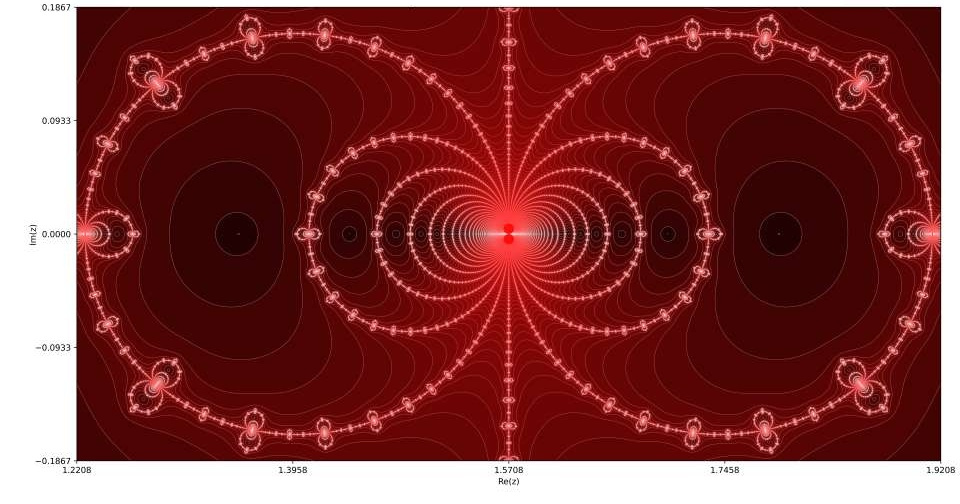}
\caption{Zoom around the pole $\pi/2$, on
$[1.2208,1.9208]\times[-0.1867,0.1867]$.}
\end{subfigure}
\caption{The Newton map $\Nsin(z)=z-\tan z$. Points are shaded according to the number of iterations needed to approach a root of $\sin z$. The darkest regions surround the superattracting fixed points $k\pi$, which lie in the immediate basins $U_k$; the bright curves accumulate on $J(\Nsin)$, which contains the vertical lines $\Re z=\pi/2+k\pi$ through the poles.}\label{fig:sin}
\end{figure}

\begin{thmx}[{\cite[Theorem~C]{BFJK25}}]\label{thm:C}
Let $f(z)=z-\tan z$, Newton's method for $g(z)=\sin z$. Then the Julia set $\mathcal{J}(f)$ is locally connected.
\end{thmx}
According to \cite{BFJK25}, this is the first non-trivial example of a
locally connected Julia set of a transcendental meromorphic map outside class $\B$ with infinitely many unbounded Fatou components. Its proof combines Theorem~\ref{thm:A}, applied to the immediate basins of $\Nsin$, with Whyburn's criterion (Lemma~\ref{lem:whyburn}); see Remark~\ref{rem:proofC}. The authors of \cite{BFJK25} also point out that most of the Newton maps of trigonometric polynomials studied in \cite{BFJK17} have infinitely many unbounded attracting basins satisfying the hypotheses of Theorem~\ref{thm:A}. The dynamical plane of $\Nsin$ is shown in Figure~\ref{fig:sin}.

\subsection{What is new in this paper}\label{subsec:new}

Since the functions $\sin z$, $\cos z$, $\sinh z$ and $\cosh z$ are so closely related, it is natural to ask the following questions.
\begin{question}\label{q:1}
Are the Julia sets of the Newton maps $\Ncos(z)=z+\cot z$, $\Nsinh(z)=z-\tanh
z$ and $\Ncosh(z)=z-\coth z$ of $\cos z$, $\sinh z$ and $\cosh z$ locally connected?
\end{question}

\begin{question}\label{q:2}
Can one check directly, and with explicit constants, that the immediate basins of $\Ncos$ satisfy the hypotheses of Theorem~\ref{thm:A}?
\end{question}
In Section~\ref{sec:main} we answer both questions affirmatively. The answer to Question~\ref{q:1} (Theorem~\ref{thm:lc}) rests on an elementary observation, Lemma~\ref{lem:affine}: Newton's method commutes with affine changes of variable. Since
\[
\cos z=\sin(z+\tfrac{\pi}{2}),\; \sinh z=-i\sin(iz)\; \mbox{and}\;
\cosh z=\sin(iz+\tfrac{\pi}{2}),
\]
the four Newton maps are affinely conjugate (Theorem~\ref{thm:conj}), and Theorem~\ref{thm:C} applies to all of them. More generally, every function $Ae^{\lambda z}+Be^{-\lambda z}$ with $A,B,\lambda\ne0$ is a constant multiple of $\sin(az+b)$ for suitable $a\ne0$ and $b$, so the same conclusion holds for the Newton maps of all these functions (Theorem~\ref{thm:general}). Concerning Question~\ref{q:2},
$\cos z$ is a trigonometric polynomial, so we do not claim the statement itself as new: it falls within the scope of the remark of \cite{BFJK25} quoted above. Our contribution (Theorem~\ref{thm:direct}) is a complete verification with explicit constants, in which the half-planes $\{\pm\Im z>2\}$ are shown to be repelling petals at infinity directly from Definition~\ref{def:petal}. Along the way we describe the dynamical planes of all four maps explicitly (Corollary~\ref{cor:structure} and Table~\ref{tab:maps}) and illustrate them numerically (Figures~\ref{fig:sin}--\ref{fig:cosh}). This note is therefore largely expository. In Remark~\ref{rem:entire} we contrast these maps with Newton's method for $e^z-1$, which is an entire function whose Julia set is not locally connected.\vspace{1.2mm}

The paper is organized as follows. Section~\ref{sec:main} states the main results and their consequences, Section~\ref{sec:lemmas} collects the auxiliary lemmas and propositions, and Section~\ref{sec:proofs} contains the proofs of the main results.

\section{Main results and their consequences}\label{sec:main}

Throughout the paper, we write
\begin{equation}\label{eq:maps}
\begin{aligned}
\Nsin(z)&=z-\tan z, &\; \Ncos(z)&=z+\cot z,\\
\Nsinh(z)&=z-\tanh z, &\; \Ncosh(z)&=z-\coth z
\end{aligned}
\end{equation}
for Newton's methods of $\sin z$, $\cos z$, $\sinh z$ and $\cosh z$; for
instance, $z-\cos z/(-\sin z)=z+\cot z$ and $z-\cosh z/\sinh z=z-\coth z$. We
also consider the affine maps
\begin{equation}\label{eq:affine}
A_1(z)=z+\frac{\pi}{2},\; A_2(z)=iz,\; A_3(z)=iz+\frac{\pi}{2},
\end{equation}
extended to $\Ch$ by $A_j(\infty)=\infty$.

\begin{theorem}\label{thm:conj}
The Newton maps \eqref{eq:maps} satisfy
\begin{equation}\label{eq:conj}
\begin{gathered}
\Ncos=A_1^{-1}\circ\Nsin\circ A_1,\;
\Nsinh=A_2^{-1}\circ\Nsin\circ A_2,\\
\Ncosh=A_3^{-1}\circ\Nsin\circ A_3=A_2^{-1}\circ\Ncos\circ A_2 .
\end{gathered}
\end{equation}
Consequently,
\[
J(\Ncos)=A_1^{-1}(J(\Nsin)),\; J(\Nsinh)=A_2^{-1}(J(\Nsin)),\;
J(\Ncosh)=A_3^{-1}(J(\Nsin)),
\]
and the maps $A_j^{-1}$ send Fatou components onto Fatou components, immediate basins onto immediate basins, critical points onto critical points of the same local degree, and poles onto poles of the same order.
\end{theorem}
More generally, the following holds.
\begin{theorem}\label{thm:general}
Let $A,B,\lambda\in\C\setminus\{0\}$ and $g(z)=Ae^{\lambda z}+Be^{-\lambda z}$. Choose $\varphi\in\C$ with $e^{2i\varphi}=-A/B$, and put \[
T(z)=-i\lambda z+\varphi,\; C=2iAe^{-i\varphi}.
\]
Then $g=C\,(\sin\circ T)$ and
\[
N_g=T^{-1}\circ\Nsin\circ T .
\]
Consequently, $J(N_g)=T^{-1}(J(\Nsin))$, the zeros of $g$ are the points
$T^{-1}(k\pi)$, $k\in\Z$, and $T^{-1}$ maps the immediate basin of $k\pi$ for $\Nsin$ onto the immediate basin of $T^{-1}(k\pi)$ for $N_g$; in particular, $N_g$ has infinitely many unbounded immediate basins of attraction.
\end{theorem}

\begin{remark}\label{rem:examples}
Theorem~\ref{thm:general} applies, for instance, to
\begin{align*}
a\sin z+b\cos z&=\frac{b-ia}{2}\,e^{iz}+\frac{b+ia}{2}\,e^{-iz}
&&(a,b\in\C,\ a^2+b^2\ne0),\\
a\cosh z+b\sinh z&=\frac{a+b}{2}\,e^{z}+\frac{a-b}{2}\,e^{-z}
&&(a,b\in\C,\ a^2\ne b^2),\\
\sin(\alpha z+\beta)&=\frac{e^{i\beta}}{2i}\,e^{i\alpha z}
-\frac{e^{-i\beta}}{2i}\,e^{-i\alpha z} &&(\alpha,\beta\in\C,\ \alpha\ne0).
\end{align*}
In particular, it covers $\cos z$, $\sinh z$ and $\cosh z$; for $g=\cos$ one can take $\varphi=\pi/2$, which gives $T=A_1$ and $C=1$. If $AB=0$, then $g$ has no zeros and $N_g$ is a translation.
\end{remark}

\begin{theorem}\label{thm:lc}
The Julia sets of $\Nsin$, $\Ncos$, $\Nsinh$ and $\Ncosh$ are locally
connected. More generally, $J(N_g)$ is locally connected for every function $g$ as in Theorem~\ref{thm:general}. In particular, every Fatou component of each of these maps has a locally connected boundary.
\end{theorem}

\begin{figure}[htbp]
\centering
\includegraphics[width=0.6\textwidth]{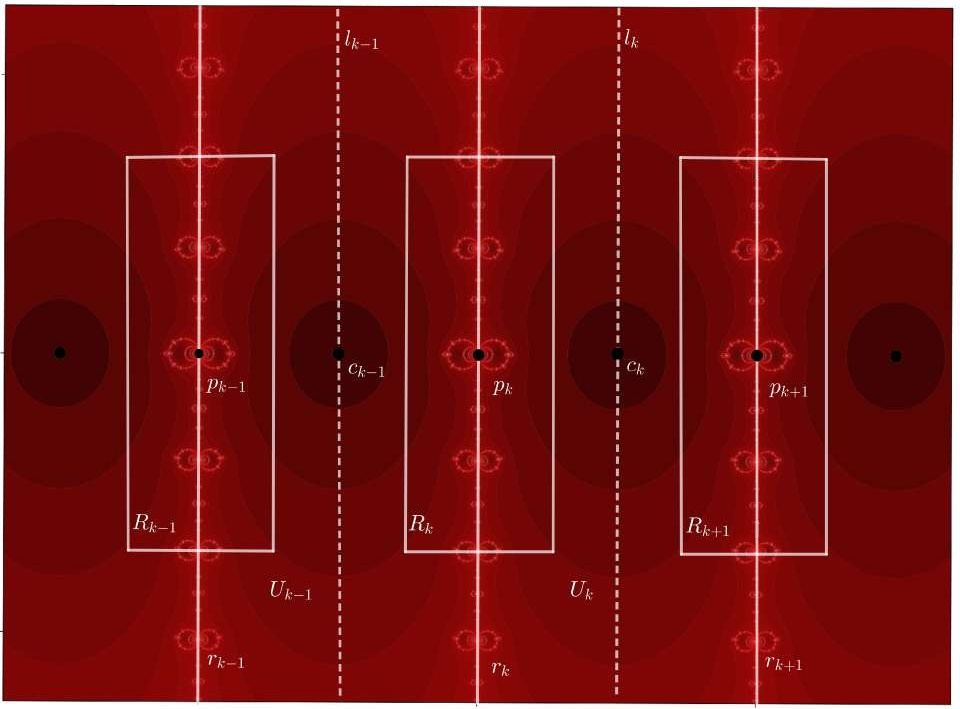}
\caption{The dynamical plane of $f=\Ncos$. The solid vertical lines
$r_k=\{\Re z=k\pi\}$ pass through the poles $p_k=k\pi$ and are contained in $\mathcal{J}(f)$; the dashed lines $\ell_k=\{\Re z=\pi/2+k\pi\}$ (labelled $l_k$ in the picture) pass through the superattracting fixed points $c_k$ and are contained in the immediate basins $U_k$. The basin $U_k$ lies between $r_k$ and $r_{k+1}$. The rectangles $R_k$ mark neighborhoods of the poles $p_k$, near which the preimages of the basins accumulate.}
\label{fig:cosannot}
\end{figure}

The next result answers Question~\ref{q:2}; see Figure~\ref{fig:cosannot} for the notation.

\begin{theorem}\label{thm:direct}
Let $f=\Ncos$, let $U_k$ be the immediate basin of the superattracting fixed point $c_k=\pi/2+k\pi$, $k\in\Z$, let $M\ge2$, and let
\[
P_+=\{z\in\C:\Im z>M\}\; \mbox{and}\; P_-=\{z\in\C:\Im z<-M\}.
\]
\begin{enumerate}[label=\textup{(\roman*)}]
\item The map $f$ is univalent on $P_\pm$, $\overline{P_\pm}\subset f(P_\pm)$, $f(P_+)\cap P_-=f(P_-)\cap P_+=\emptyset$, and $P_+$ and $P_-$ are repelling petals of $f$ at infinity in the sense of Definition~\ref{def:petal}, with $g\equiv1$.\vspace{1.2mm}

\item For every $k\in\Z$, the basin $U_k$ satisfies hypotheses \textup{(a)}, \textup{(b)} and \textup{(c)} of Theorem~\ref{thm:A}. More precisely, the set in \textup{(a)} is empty, \textup{(b)} holds with $L=\{c_k\}$, and $\overline{U_k}\setminus D\subset P_+\cup P_-$ for the disc $D$ centred at $0$ of radius $(|k|+1)\pi+M+1$.
\end{enumerate}
Consequently, $\partial U_k$ is locally connected for every $k\in\Z$, and the conclusion of Theorem~\ref{thm:B} holds for $U_k$.
\end{theorem}

\begin{remark}\label{rem:L}
The compact set $L$ in Theorem~\ref{thm:A}(b) cannot be omitted here. Indeed, $c_k\in f^{-n}(c_k)$ and $\deg_{c_k}f^n=3^n\to\infty$, so the supremum in Theorem~\ref{thm:A}(b) is infinite for $z=c_k$; the hypothesis holds only because $c_k\in L$.
\end{remark}

\begin{table}[htbp]
\centering
\renewcommand{\arraystretch}{1.25}
\small
\setlength{\tabcolsep}{4pt}
\begin{tabular}{@{}lllll@{}}
\toprule
 & $\Nsin$ & $\Ncos$ & $\Nsinh$ & $\Ncosh$\\
\midrule
$N(z)$ & $z-\tan z$ & $z+\cot z$ & $z-\tanh z$ & $z-\coth z$\\
fixed points $c_k$ & $k\pi$ & $\frac{\pi}{2}+k\pi$ & $k\pi i$ &
 $(k+\frac12)\pi i$\\
poles $p_k$ & $\frac{\pi}{2}+k\pi$ & $k\pi$ & $(k+\frac12)\pi i$ & $k\pi i$\\
lines $\ell_k\subset U_k$ & $\Re z=k\pi$ & $\Re z=\frac{\pi}{2}+k\pi$ &
 $\Im z=k\pi$ & $\Im z=(k+\frac12)\pi$\\
lines $r_k\subset J$ & $\Re z=\frac{\pi}{2}+k\pi$ & $\Re z=k\pi$ &
 $\Im z=(k+\frac12)\pi$ & $\Im z=k\pi$\\
period $\tau$ & $\pi$ & $\pi$ & $\pi i$ & $\pi i$\\
petals $P_\pm$ & $\pm\Im z>M$ & $\pm\Im z>M$ & $\pm\Re z>M$ & $\pm\Re z>M$\\
conjugacy & --- & $A_1$ & $A_2$ & $A_3$\\
\bottomrule
\end{tabular}
\caption{The four Newton maps ($k\in\Z$, $M\ge2$). The last row gives the map $A_j$ in \eqref{eq:affine} with $N=A_j^{-1}\circ\Nsin\circ A_j$.}\label{tab:maps}
\end{table}

\begin{corollary}\label{cor:structure}
Let $f$ be one of the four maps \eqref{eq:maps}, and let $c_k$, $p_k$,
$\ell_k$, $r_k$, $\tau$ and $P_\pm$ be as in Table~\ref{tab:maps}. Then the following hold.
\begin{enumerate}
\item The critical points of $f$ are exactly the points $c_k$, $k\in\Z$; each $c_k$ is a superattracting fixed point with $\deg_{c_k}f=3$. The map $f$ has no
finite asymptotic values, so that $\Sing(f)=\PS(f)=\{c_k:k\in\Z\}$. The poles
of $f$ are the points $p_k$, and they are simple.
\item $f(z+\tau)=f(z)+\tau$ for every $z\in\C$.
\item For every $k\in\Z$, the immediate basin $U_k$ of $c_k$ is simply
connected, $U_k=U_0+k\tau$, $\deg(f|_{U_k})=3$, and $U_k$ contains the line $\ell_k$, which satisfies $f(\ell_k)\subset\ell_k$. The lines $r_j$, $j\in\Z$, are contained in $\mathcal{J}(f)$ and satisfy $f(r_j\setminus\{p_j\})\subset r_j$. The
basin $U_k$ lies in the open strip bounded by the two lines $r_j$ adjacent to $\ell_k$; its boundary contains exactly two poles of $f$, namely the poles on these two lines, and both are accessible from $U_k$. Moreover, $U_k$ has exactly two accesses to $\infty$.
\item $\mathcal{J}(f)\cap\C=\bigcup_{j\in\Z}J_j$, where $J_j$ is the connected component of $\mathcal{J}(f)\cap\C$ containing $r_j$.
\item Every Fatou component of $f$ is eventually mapped into some $U_k$; in particular, $f$ has no wandering domains.
\item For every $M\ge2$, the half-planes $P_+$ and $P_-$ are repelling petals of $f$ at infinity, and every $U_k$ is tame at infinity with respect to them. 
\item $\mathcal{J}(f)$ is locally connected, and the boundary of every Fatou component of $f$ is locally connected.
\end{enumerate}
\end{corollary}

\begin{remark}\label{rem:classB}
For each of the four maps, the set $\Sing(f)=\{c_k:k\in\Z\}$ is unbounded, so $f\notin\B$. Moreover, $\PS(f)=\{c_k:k\in\Z\}$ is not compact, since it accumulates at $\infty\in \mathcal{J}(f)$. Hence $f$ is not hyperbolic in the usual sense that $\overline{\PS(f)}$ is a compact subset of $\mathcal{F}(f)$ (cf.\ \cite{BFR15}), and it is not geometrically finite, because $\Sing(f)\cap \mathcal{F}(f)=\Sing(f)$ is not compact. Thus the results for maps in class $\B$ and for geometrically finite maps recalled in Section~\ref{subsec:background} do not apply, and Theorem~\ref{thm:A} is needed. Note that the functions $\sin z$, $\cos z$, $\sinh z$ and $\cosh z$ themselves
belong to $\B$ (their singular values are $\pm1$ or $\pm i$), but this plays no role: what matters is the class of the Newton map.
\end{remark}

\begin{remark}\label{rem:entire}
The results above should be contrasted with the case of Newton maps that are
entire. For instance, Newton's method for $g(z)=e^z-1$ is the transcendental
entire function
\[
N_g(z)=z-\frac{e^z-1}{e^z}=z-1+e^{-z}.
\]
Its fixed points are the zeros $2\pi ik$, $k\in\Z$, of $g$, and they are
superattracting, since $N_g'(z)=1-e^{-z}$. For real $x>0$ we have
$1-x<e^{-x}<1$, hence $0<N_g(x)<x$; thus the orbit of $x$ decreases to a fixed
point of $N_g$ in $[0,x)$, that is, to $0$. Hence $[0,+\infty)$ is contained in
the immediate basin $U$ of $0$, which is therefore unbounded. Since $U$ is an
invariant attracting basin and not a Baker domain, its boundary is not locally
connected by \cite{BW91,Kis97}, and $J(N_g)$ is not locally connected by
\cite{BD00,Osb13}; see Section~\ref{subsec:background}.
\end{remark}

\begin{remark}\label{rem:novelty} Theorem~\ref{thm:lc} follows from Theorem~\ref{thm:C} by elementary conjugacies (Theorems~\ref{thm:conj} and~\ref{thm:general}), so its essential input is the work of \cite{BFJK25}. Theorem~\ref{thm:direct}, on the other hand, does not use Theorem~\ref{thm:C}: it shows independently that every immediate basin of $\Ncos$ has a locally connected boundary. To deduce the local connectivity of the whole Julia set $J(\Ncos)$ without using the conjugacy, one would still have to verify both conditions of Whyburn's criterion (Lemma~\ref{lem:whyburn}) for all Fatou components, along the lines of \cite[Section~8]{BFJK25}; see Remark~\ref{rem:proofC}.
\end{remark}

\begin{figure}[htbp]
\centering
\begin{subfigure}{0.6\textwidth}
\centering
\includegraphics[width=\textwidth]{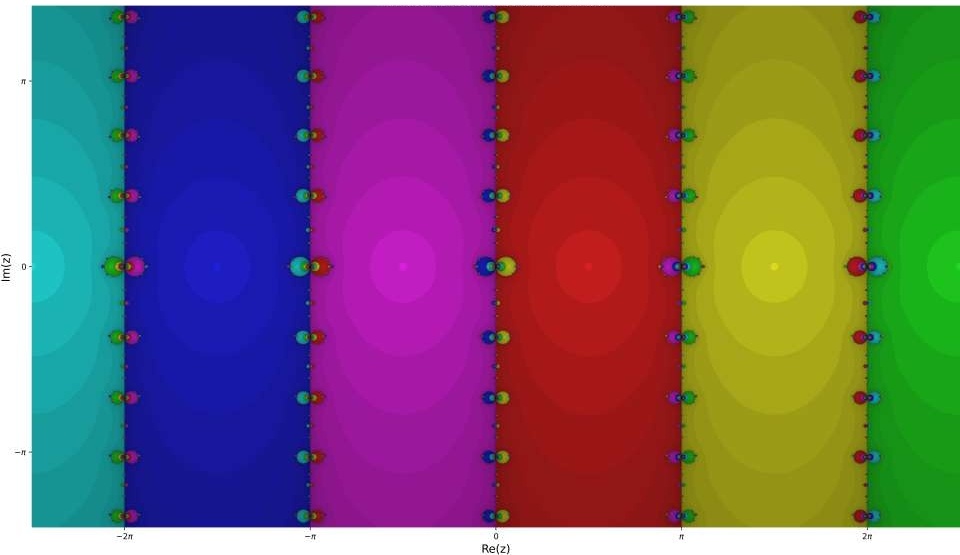}
\caption{Basins of the roots $\pi/2+k\pi$ of $\cos z$.}
\end{subfigure}

\begin{subfigure}{0.6\textwidth}
\centering
\includegraphics[width=\textwidth]{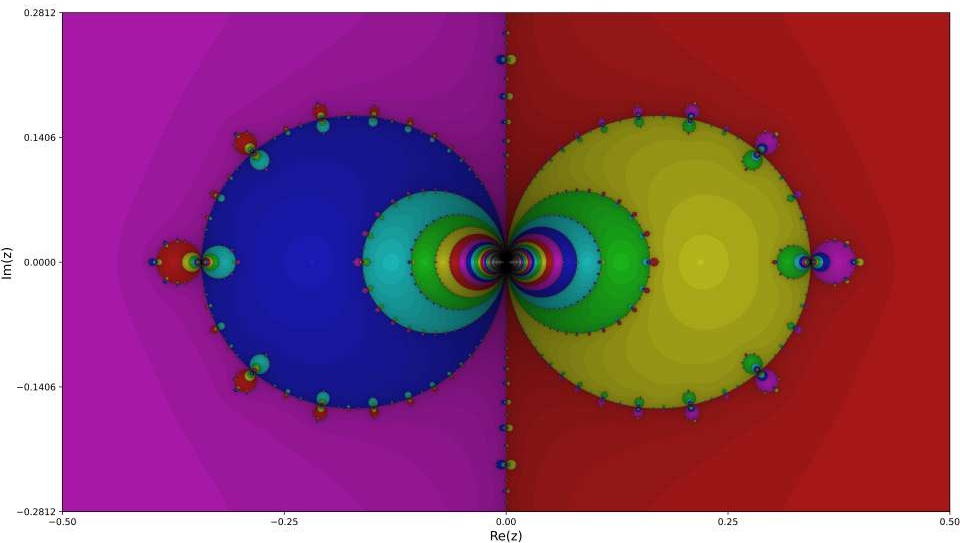}
\caption{Zoom around the pole $p_0=0$, on $[-0.5,0.5]\times[-0.2812,0.2812]$.}
\end{subfigure}
\caption{Each point is coloured according to the root of $\cos z$ to which the iteration of $\Ncos$ converges. The coloured regions are the basins of attraction, and their common boundary is the Julia set $J(\Ncos)$.} \label{fig:cosbasins}
\end{figure}

\begin{figure}[htbp]
\centering
\begin{subfigure}{0.6\textwidth}
\centering
\includegraphics[width=\textwidth]{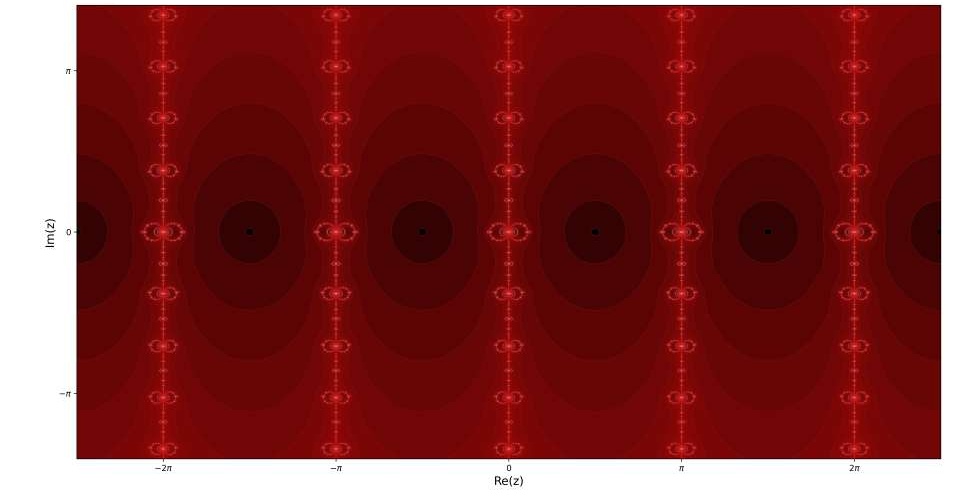}
\caption{The dynamical plane of $\Ncos$.}
\end{subfigure}

\medskip
\begin{subfigure}{0.6\textwidth}
\centering
\includegraphics[width=\textwidth]{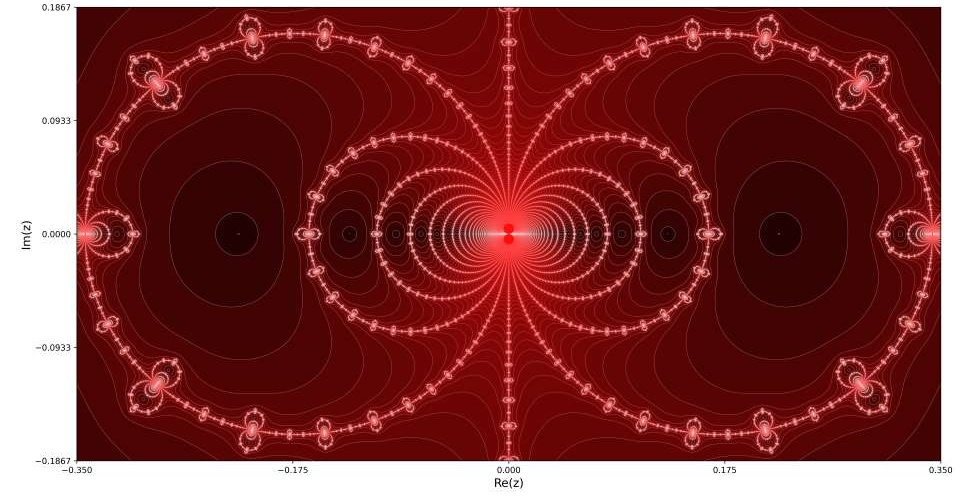}
\caption{Zoom around the pole $p_0=0$, on
$[-0.35,0.35]\times[-0.1867,0.1867]$.}
\end{subfigure}
\caption{The Newton map $\Ncos(z)=z+\cot z$, shaded according to the number of
iterations needed to approach a root of $\cos z$; compare
Figure~\ref{fig:sin}.}
\label{fig:cos}
\end{figure}

\begin{figure}[htbp]
\centering
\begin{subfigure}{0.48\textwidth}
\centering
\includegraphics[width=\textwidth]{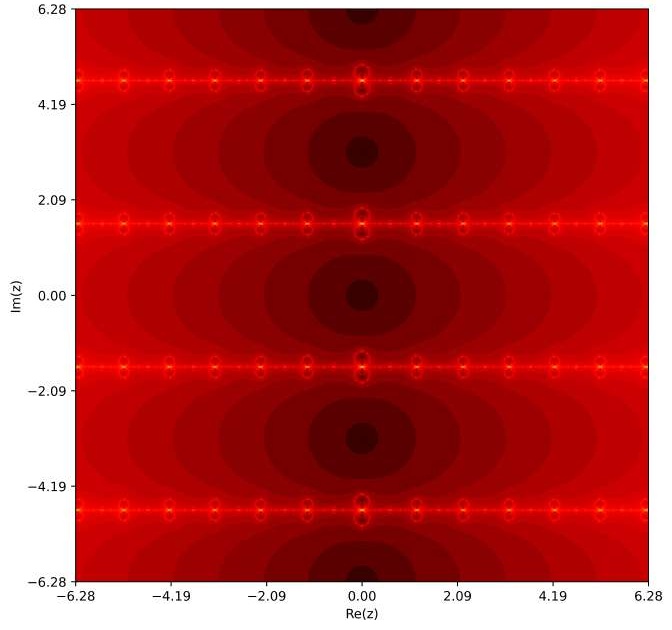}
\caption{The dynamical plane of $\Nsinh$.}
\end{subfigure}\hfill
\begin{subfigure}{0.48\textwidth}
\centering
\includegraphics[width=\textwidth]{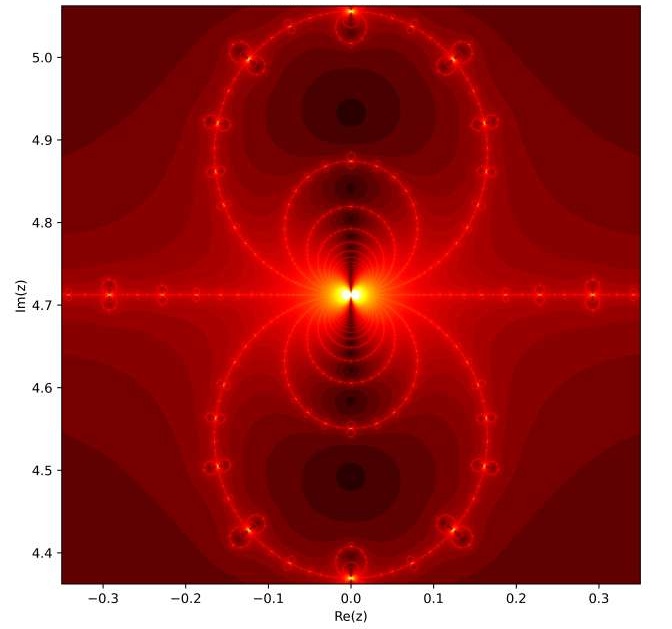}
\caption{Zoom around the pole $p_1=3\pi i/2\approx4.71i$.}
\end{subfigure}
\caption{The Newton map $\Nsinh(z)=z-\tanh z$. Its horizontal strips are
obtained from the vertical strips of $\Nsin$ by the rotation $A_2^{-1}$; the
lines $\Im z=(k+\frac12)\pi$ are contained in $J(\Nsinh)$.}
\label{fig:sinh}
\end{figure}

\begin{figure}[htbp]
\centering
\begin{subfigure}{0.48\textwidth}
\centering
\includegraphics[width=\textwidth]{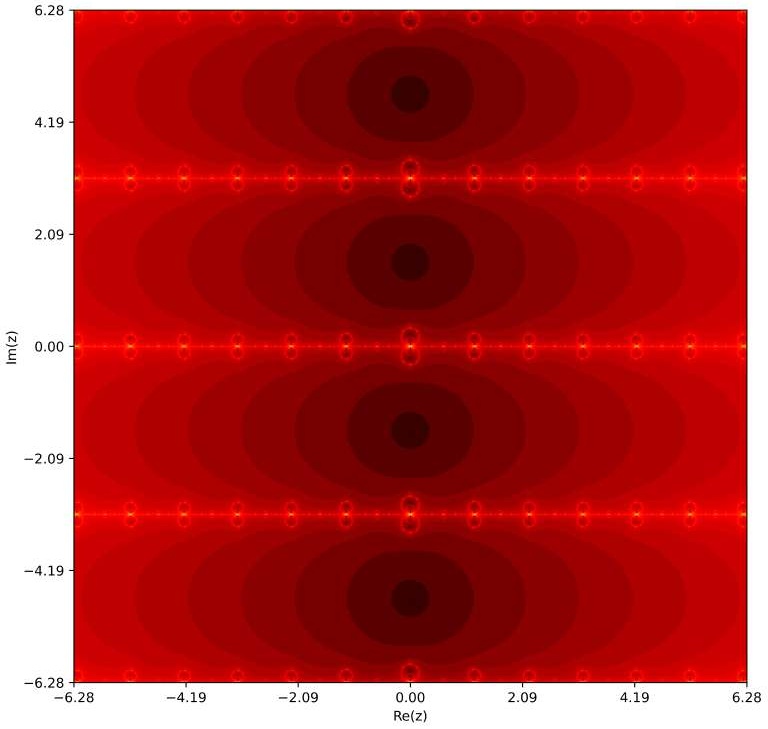}
\caption{The dynamical plane of $\Ncosh$.}
\end{subfigure}\hfill
\begin{subfigure}{0.48\textwidth}
\centering
\includegraphics[width=\textwidth]{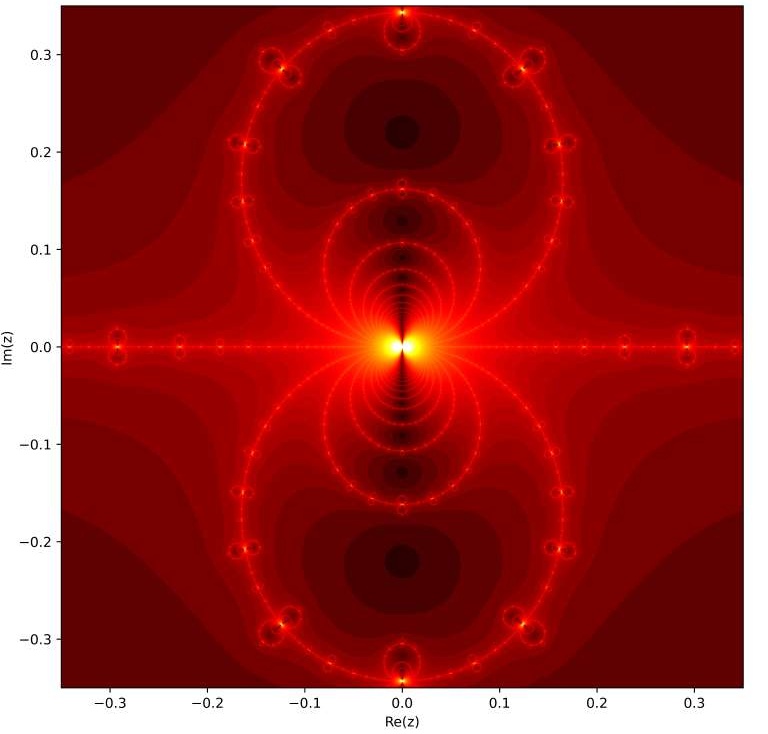}
\caption{Zoom around the pole $p_0=0$.}
\end{subfigure}
\caption{The Newton map $\Ncosh(z)=z-\coth z$. Its horizontal strips are
obtained from the vertical strips of $\Ncos$ by the rotation $A_2^{-1}$; the
lines $\Im z=k\pi$ are contained in $J(\Ncosh)$.}
\label{fig:cosh}
\end{figure}

\begin{remark}\label{rem:figures}
Figures~\ref{fig:sin}--\ref{fig:cosh} were produced in Python by iterating the Newton maps on a grid of initial points. In Figure~\ref{fig:cosbasins} each point is coloured according to the root to which it converges; in the other figures it is shaded according to the number of iterations needed to approach a root, so that slowly converging points, which accumulate on the Julia set, appear bright. The pictures for $\Ncos$ (Figures~\ref{fig:cosannot}, \ref{fig:cosbasins} and~\ref{fig:cos}) are those for $\Nsin$ translated by $\pi/2$, and the conjugacy $A_2(z)=iz$ rotates the vertical strips of $\Nsin$ and $\Ncos$ by $90^\circ$ into the horizontal strips of $\Nsinh$ and $\Ncosh$ (Figures~\ref{fig:sinh} and~\ref{fig:cosh}).
\end{remark}

\section{Auxiliary lemmas and propositions}\label{sec:lemmas}

\subsection{Topological tools}

\begin{lemma}[Torhorst's theorem, {\cite[Theorem~2.2]{Why63}}]
\label{lem:torhorst} Let $J\subset\Ch$ be a compact, connected and locally connected set. Then the boundary of every connected component of $\Ch\setminus J$ is locally connected.
\end{lemma}
\begin{lemma}[Whyburn's criterion, {\cite[Theorem~4.4,~p.~113]{Why63}}]
\label{lem:whyburn} A compact connected set $J\subset\Ch$ is locally connected if, and only if, the following two conditions hold:
\begin{enumerate}
\item the boundary of every connected component of $\Ch\setminus J$ is locally connected; \item for every $\varepsilon>0$, only finitely many connected components of $\Ch\setminus J$ have spherical diameter greater than $\varepsilon$.
\end{enumerate}
\end{lemma}

\begin{remark}\label{rem:proofC}
In \cite[Section~8]{BFJK25}, Theorem~\ref{thm:C} is deduced from
Lemma~\ref{lem:whyburn} applied to $J=J(\Nsin)$, which is connected by
\cite{BFJK14}. Condition (a) is obtained by applying Theorem~\ref{thm:A} to
the immediate basins $U_k$ and then treating the remaining Fatou components,
each of which is a preimage of some $U_k$ (Lemma~\ref{lem:sin}(e)). For
condition (b), the main step is to show that
$\sum_U\diam_{\mathrm{sph}}(U)^2<\infty$, the sum being taken over all Fatou
components $U$. Roughly speaking, one proves that the branches of $f^{-n}$ on
the Fatou components contained in $f^{-1}(U_k)\setminus U_k$, $k\in\Z$, have
uniformly bounded spherical distortion, so that the ratio between the square of
the spherical diameter and the spherical area of a Fatou component is
controlled uniformly along its iterated preimages; since distinct Fatou
components are disjoint, their total spherical area is at most the area $4\pi$
of $\Ch$. The summability of the squared diameters implies that, for every
$\varepsilon>0$, only finitely many Fatou components have spherical diameter
at least $\varepsilon$. In the meromorphic setting this step is delicate,
because preperiodic Fatou components may accumulate at poles and prepoles; we
refer to \cite[Section~8]{BFJK25} for the details.
\end{remark}

\subsection{Newton's method and affine changes of variable}

\begin{lemma}[Affine invariance of Newton's method]\label{lem:affine}
Let $g$ and $h$ be non-constant entire functions, let $A(z)=az+b$ with
$a,b\in\C$, $a\ne0$, and let $c\in\C\setminus\{0\}$. If $g\circ A=c\,h$, then
\[
N_h=A^{-1}\circ N_g\circ A .
\]
\end{lemma}

\begin{proof}
Differentiating $g(A(z))=c\,h(z)$ gives $a\,g'(A(z))=c\,h'(z)$. Hence, for
every $z$ with $h'(z)\ne0$,
\begin{align*}
	A^{-1}\bigl(N_g(A(z))\bigr)
	&=\frac{1}{a}\Bigl(az+b-\frac{g(A(z))}{g'(A(z))}-b\Bigr)\\&
	=z-\frac{g(A(z))}{a\,g'(A(z))}\\&
	=z-\frac{c\,h(z)}{c\,h'(z)}\\& =N_h(z).
\end{align*}
Both sides are meromorphic in $\C$ and agree outside a discrete set, so they
agree everywhere.
\end{proof}

\begin{lemma}[Affine conjugacy]\label{lem:conjugacy}
	Let $f\colon\C\to\Ch$ be a nonconstant meromorphic function. Let
	$A(z)=az+b$, where $a\ne0$, extend $A$ by $A(\infty)=\infty$, and
	set $F=A^{-1}\circ f\circ A$. Then the following statements hold.
	\begin{enumerate}
		\item For every $n\geq1$, the iterate $F^n$ is defined at $z$ if
		and only if $f^n$ is defined at $A(z)$. Wherever these iterates
		are defined, we have
		\begin{equation}\label{eq:affine-iterates}
			F^n(z)=A^{-1}\bigl(f^n(A(z))\bigr).
		\end{equation}
		
		\item The Fatou and Julia sets satisfy
		\begin{equation}\label{eq:affine-fatou-julia}
			\mathcal{F}(F)=A^{-1}\bigl(\mathcal{F}(f)\bigr),
			\qquad
			\mathcal{J}(F)=A^{-1}\bigl(\mathcal{J}(f)\bigr).
		\end{equation}
		The map $A^{-1}$ sends each Fatou component of $f$ onto a Fatou
		component of $F$. If $c\in\C$ is an attracting fixed point of
		$f$, then $A^{-1}(c)$ is an attracting fixed point of $F$ with the same multiplier, and $A^{-1}$ maps their immediate basins onto one another.\vspace{1.2mm}
		
		\item For each $z\in\C$, we have
		$\deg_z F=\deg_{A(z)}f$. In particular, $z$ is a critical point
		of $F$ if and only if $A(z)$ is a critical point of $f$. A point
		$w\in\C$ is a pole of $F$ of order $m$ if and only if $A(w)$ is a pole of $f$ of order $m$. A finite value $v$ is asymptotic for $F$	if and only if $A(v)$ is asymptotic for $f$.\vspace{1.2mm}
		
		\item The Julia set $\mathcal{J}(F)$ is locally connected if and
		only if $\mathcal{J}(f)$ is locally connected. If $U$ is a Fatou
		component of $f$, then $\partial A^{-1}(U)$ is locally connected
		if and only if $\partial U$ is locally connected. Moreover,
		$A^{-1}$ induces a bijection between the accesses to $\infty$
		from $U$ and those from $A^{-1}(U)$.
	\end{enumerate}
\end{lemma}

\begin{proof}
	{(a).}	For $n=1$, \eqref{eq:affine-iterates} is the definition of $F$. Suppose that the identity and the correspondence of domains hold for some $n\geq1$. Since $A$ fixes $\infty$, the two values
	$F^n(z)$ and $f^n(A(z))$ are finite simultaneously. If $f$ is
	transcendental, an orbit that reaches $\infty$ cannot be iterated
	further. If $f$ is rational, we use its extension to $\Ch$, and
	the conjugacy also holds there. Thus the next iterate is defined
	on one side exactly when it is defined on the other. At every
	such point, we obtain
	\[
	F^{n+1}(z)
	=A^{-1}\!\left(f\bigl(A(F^n(z))\bigr)\right)
	=A^{-1}\!\left(f\bigl(f^n(A(z))\bigr)\right).
	\]
	This proves \eqref{eq:affine-iterates} by induction.\vspace{1.2mm}
	
	{(b).}
	Let $V$ be a neighborhood of $z$ on which all iterates of $F$
	are defined. By part~(a), all iterates of $f$ are defined on
	$A(V)$. If a subsequence $f^{n_j}$ converges locally uniformly
	in the spherical metric on $A(V)$ to a limit $g$, then
	\eqref{eq:affine-iterates} gives
	\[
	F^{n_j}\longrightarrow A^{-1}\circ g\circ A
	\]
	locally uniformly in the spherical metric on $V$. Indeed, $A$
	maps compact subsets of $V$ to compact subsets of $A(V)$, and
	$A^{-1}$ is uniformly continuous on $\Ch$. Applying the same
	argument with $A^{-1}$ in place of $A$ proves the converse.
	Hence normality is equivalent on the two neighborhoods, which
	gives \eqref{eq:affine-fatou-julia}. Since $A^{-1}$ is a
	homeomorphism, it maps connected components onto connected
	components.
	
	Now let $f(c)=c$ and put $\widetilde c=A^{-1}(c)$. We have
	$F(\widetilde c)=\widetilde c$, and the chain rule gives
	\begin{equation}\label{eq:affine-multiplier}
		F'(\widetilde c)
		=(A^{-1})'(c)\,f'(c)\,A'(\widetilde c)
		=\frac1a f'(c)a=f'(c).
	\end{equation}
	Thus $\widetilde c$ is attracting whenever $c$ is. By
	\eqref{eq:affine-iterates}, an orbit converges to $\widetilde c$
	under $F$ exactly when its image under $A$ converges to $c$
	under $f$. Taking the components of the two attracting basins
	that contain their fixed points gives the immediate-basin
	correspondence.\vspace{1.2mm}
	
	{(c).}
	In local coordinates centered at a point and its image, the
	local degree is the order of the first nonzero term of the
	map. Pre- and post-composition by biholomorphic coordinate
	changes preserve this order. Since $A$ and $A^{-1}$ are
	biholomorphic on $\Ch$, we obtain $\deg_z F=\deg_{A(z)}f$,
	including at points whose image is $\infty$, where we use
	a coordinate centered at $\infty$.
	
	Because $A^{-1}(\infty)=\infty$, we have $F(w)=\infty$ exactly
	when $f(A(w))=\infty$. At a pole, the local degree equals the
	pole order. The local-degree identity therefore gives the
	equality of pole orders.
	
	If $\gamma(t)\to\infty$ and $F(\gamma(t))\to v\in\C$, then
	$A(\gamma(t))\to\infty$ and
	\[
	f(A(\gamma(t)))=A(F(\gamma(t)))\longrightarrow A(v).
	\]
	Thus $A(v)$ is an asymptotic value of $f$. Conversely, applying
	$A^{-1}$ to a curve tending to infinity along which $f$ tends
	to $A(v)$ gives an asymptotic curve for $F$ with value $v$.\vspace{1.2mm}
	
	{(d).}
	By \eqref{eq:affine-fatou-julia}, the restriction of $A^{-1}$
	is a homeomorphism from $\mathcal{J}(f)$ onto
	$\mathcal{J}(F)$. For any Fatou component $U$, it also maps
	$\partial U$ homeomorphically onto $\partial A^{-1}(U)$,
	where we take boundaries in $\Ch$. Local connectivity is
	preserved under these homeomorphisms.\vspace{1.2mm}
	
	Finally, fix $z_0\in U$. An access to infinity is represented
	by a curve from $z_0$ to $\infty$ whose points before the
	endpoint lie in $U$, up to homotopy through such curves with
	the endpoints fixed. Composing the curve with $A^{-1}$ gives
	a curve from $A^{-1}(z_0)$ to $\infty$ in the corresponding
	component. Composing a homotopy with $A^{-1}$ preserves the
	endpoint and interior conditions, so the correspondence is
	well defined on access classes. Composition with $A$ gives
	the inverse correspondence.
\end{proof}
\subsection{The Newton map of \texorpdfstring{$\sin z$}{sin z}}

We shall use the following properties of $\Nsin$, which are proved in
\cite[Example~7.2]{BFJK17}, \cite[Proposition~4.1]{BFJK20} and
\cite[Section~8]{BFJK25}.

\begin{lemma}\label{lem:sin}
Let $f=\Nsin$. Then the following hold.
\begin{enumerate}
\item $f(z+\pi)=f(z)+\pi$ for $z\in\C$. The map $f$ has infinitely many simply
connected immediate basins of attraction $U_k$, $k\in\Z$, of the
superattracting fixed points $c_k=k\pi$, such that $U_k=U_0+k\pi$ and
$\deg_{c_k}f=\deg(f|_{U_k})=3$. The points $c_k$ are the only critical points
of $f$, and the vertical lines $\ell_k=\{k\pi+it:t\in\R\}$ are $f$-invariant and
contained in $U_k$. Moreover, $f$ has no finite asymptotic values.
\item The poles $p_k=\pi/2+k\pi$, $k\in\Z$, of $f$ are simple, and they lie on
the vertical lines $r_k=\{\pi/2+k\pi+it:t\in\R\}$, which are contained in
$\mathcal{J}(f)$.
\item As $\Im z\to\pm\infty$,
\[
f(z)=\begin{cases}
z-i+O\bigl(e^{-2\Im z}\bigr), & \Im z\to+\infty,\\
z+i+O\bigl(e^{2\Im z}\bigr), & \Im z\to-\infty.
\end{cases}
\]
\item $\mathcal{J}(f)\cap\C=\bigcup_{k\in\Z}J_k$, where $J_k$ is the connected component
of $\mathcal{J}(f)\cap\C$ containing $r_k$. For every $k\in\Z$, the basin $U_k$ has
exactly two accesses to infinity, and $\partial U_k$ contains exactly two poles
of $f$, namely $p_{k-1}$ and $p_k$, which are accessible from $U_k$.
\item All Fatou components of $f$ are preperiodic and are eventually mapped
into some $U_k$. In particular, $f$ has no wandering domains.
\end{enumerate}
\end{lemma}
The elementary parts of Lemma~\ref{lem:sin} are checked exactly as for
$\Ncos$ in Lemma~\ref{lem:cos} below: $\tan$ is $\pi$-periodic, $f'(z)=1-\sec^2
z=-\tan^2 z$ vanishes exactly at the points $k\pi$,
$f(k\pi+it)=k\pi+i(t-\tanh t)$, where $|t-\tanh t|<|t|$ for $t\ne0$, and, since $\tan(\pi/2+k\pi+it)=i\coth t$, we have
$f(\pi/2+k\pi+it)=\pi/2+k\pi+i(t-\coth t)$ for $t\ne0$; in particular,
$f(r_k\setminus\{p_k\})\subset r_k$.

\begin{remark}\label{rem:realtrace}
Cloutier \emph{et al.} \cite{CDBSW24} give an elementary
study of $\Nsin$, including its symmetries and its behaviour on the real line.
In particular, $U_0\cap\R=(-x^*,x^*)$, where $x^*=1.1655\ldots$ is the solution
of $\tan x=2x$ in $(0,\pi/2)$. The points $\pm x^*$ form a repelling cycle of
period $2$ of $\Nsin$: indeed, $\Nsin(\pm x^*)=\pm(x^*-\tan x^*)=\mp x^*$, and
the multiplier of this cycle is $\Nsin'(x^*)\,\Nsin'(-x^*)=\tan^4x^*
=16(x^*)^4>1$.
\end{remark}

\subsection{Elementary estimates}

For $t>0$ we put
\begin{equation}\label{eq:eta}
\eta(t)=\frac{2}{e^{2t}-1}=\frac{e^{-t}}{\sinh t}.
\end{equation}

\begin{lemma}\label{lem:estimates}
Let $z=x+iy$ with $x,y\in\R$, $y\ne0$. Then:
\begin{enumerate}[label=\textup{(\roman*)}]
\item $|\sin z|^2=\sin^2x+\sinh^2y$ and $|\cos z|^2=\cos^2x+\sinh^2y$;
\item $|\csc^2 z|\le\sinh^{-2}y$ and $|\sec^2 z|\le\sinh^{-2}y$;
\item if $y>0$, then $|\cot z+i|\le\eta(y)$ and $|\tan z-i|\le\eta(y)$; if
$y<0$, then $|\cot z-i|\le\eta(-y)$ and $|\tan z+i|\le\eta(-y)$;
\item $\eta$ is decreasing, $\eta(t)\le4e^{-2t}$ for $t\ge1$, and
\[
\eta(2)<0.04,\; \eta(5/2)<0.014,\;\sinh^{-2}(3/2)<0.23,\; \mbox{and}\;
\sinh^{-2}(5/2)<0.03 .
\]
\end{enumerate}
\end{lemma}

\begin{proof}
	We recall that $\eta(t)=2/(e^{2t}-1)$ for $t>0$.
	
	{(i).}
	By the addition formulas, we have
	\begin{align*}
		\sin(x+iy)&=\sin x\cosh y+i\cos x\sinh y,\\
		\cos(x+iy)&=\cos x\cosh y-i\sin x\sinh y.
	\end{align*}
	Taking squared moduli and using $\cosh^2y=1+\sinh^2y$, we obtain
	\begin{align}
		|\sin z|^2
		&=\sin^2x\cosh^2y+\cos^2x\sinh^2y
		=\sin^2x+\sinh^2y,\label{eq:sin-modulus}\\
		|\cos z|^2
		&=\cos^2x\cosh^2y+\sin^2x\sinh^2y
		=\cos^2x+\sinh^2y.\label{eq:cos-modulus}
	\end{align}
	
	{(ii).}
	Since $y\ne0$, we have $\sinh^2y>0$. By
	\eqref{eq:sin-modulus} and \eqref{eq:cos-modulus}, both
	$|\sin z|^2$ and $|\cos z|^2$ are at least $\sinh^2y$.
	Taking reciprocals gives
	\[
	|\csc^2z|=\frac{1}{|\sin z|^2}\leq\sinh^{-2}y\; \mbox{and}\;
	|\sec^2z|=\frac{1}{|\cos z|^2}\leq\sinh^{-2}y.
	\]
	
	{(iii).}
	We put $q=e^{2iz}$, so $|q|=e^{-2y}$. The exponential
	representations of sine and cosine give
	\[
	\cot z=i\frac{e^{iz}+e^{-iz}}{e^{iz}-e^{-iz}}
	=i\frac{q+1}{q-1}.
	\]
	Consequently, we have
	\begin{equation}\label{eq:cot-shift-identities}
		\cot z+i=\frac{2iq}{q-1}\; \mbox{and}\;
		\cot z-i=\frac{2i}{q-1}.
	\end{equation}
	If $y>0$, then $|q|<1$. By the reverse triangle inequality,
	$|q-1|\geq1-|q|$. The first identity in
	\eqref{eq:cot-shift-identities} therefore gives
	\[
	|\cot z+i|
	\leq\frac{2|q|}{1-|q|}
	=\frac{2e^{-2y}}{1-e^{-2y}}
	=\eta(y).
	\]
	If $y<0$, then $|q|>1$, and we instead use
	$|q-1|\geq|q|-1$. The second identity in
	\eqref{eq:cot-shift-identities} gives
	\[
	|\cot z-i|
	\leq\frac{2}{|q|-1}
	=\frac{2}{e^{-2y}-1}
	=\eta(-y).
	\]
	Finally, we use $\tan z=-\cot(z+\pi/2)$. Since the shift by
	$\pi/2$ leaves the imaginary part unchanged, we apply the
	cotangent estimates at $z+\pi/2$ to obtain
	\[
	|\tan z-i|\leq\eta(y),\; y>0\; \mbox{and}\;
	|\tan z+i|\leq\eta(-y),\;y<0.
	\]
	
	{(iv).}
	Differentiating $\eta$, we find
	\[
	\eta'(t)=-\frac{4e^{2t}}{(e^{2t}-1)^2}<0,\;t>0.
	\]
	For $t\geq1$, we have $e^{-2t}\leq e^{-2}<1/2$, and hence
	\begin{equation}\label{eq:eta-decay}
		\eta(t)=\frac{2e^{-2t}}{1-e^{-2t}}\leq4e^{-2t}.
	\end{equation}
	
	To verify the numerical bounds with strict inequalities, we use
	$e^3>20$, $e^4>51$, and $e^5>144$. These follow by truncating
	the exponential series after the terms of degrees $8$, $7$,
	and $10$, respectively. Thus
	\[
	\eta(2)=\frac{2}{e^4-1}<\frac{2}{50}=0.04\; \mbox{and}\;
	\eta(5/2)=\frac{2}{e^5-1}<\frac{2}{143}<0.014.
	\]
	For the remaining bounds, we write
	\[
	\sinh^{-2}t=\frac{4e^{2t}}{(e^{2t}-1)^2}.
	\]
	The function $s\mapsto4s/(s-1)^2$ is strictly decreasing for
	$s>1$, since its derivative is $-4(s+1)/(s-1)^3$. Substituting
	the lower bounds for $e^3$ and $e^5$, we obtain
	\[
	\sinh^{-2}(3/2)<\frac{80}{361}<0.23\; \mbox{and}\;
	\sinh^{-2}(5/2)<\frac{576}{20449}<0.03.
	\]
	This completes the proof.
\end{proof}
\subsection{The Newton map of \texorpdfstring{$\cos z$}{cos z}}

\begin{lemma}\label{lem:cos}
Let $f=\Ncos$ and, for $k\in\Z$, let $c_k=\pi/2+k\pi$, $p_k=k\pi$,
$\ell_k=\{\pi/2+k\pi+it:t\in\R\}$ and $r_k=\{k\pi+it:t\in\R\}$. Then the
following hold.
\begin{enumerate}
\item $f(z+\pi)=f(z)+\pi$ and $f(\overline{z})=\overline{f(z)}$ for $z\in\C$.
\item $f'(z)=1-\csc^2z=-\cot^2z$. The critical points of $f$ are exactly the
points $c_k$, $k\in\Z$. Each $c_k$ is a superattracting fixed point, and
\[
f(c_k+w)=c_k-\frac{w^3}{3}+O(w^5)\qquad(w\to0),
\]
so that $\deg_{c_k}f=3$.
\item The poles of $f$ are the points $p_k$, and they are simple. For
$t\in\R$,
\[
f\Bigl(\frac{\pi}{2}+k\pi+it\Bigr)=\frac{\pi}{2}+k\pi+i(t-\tanh t),
\qquad f(k\pi+it)=k\pi+i(t-\coth t)\quad(t\ne0).
\]
In particular, $f(\ell_k)\subset\ell_k$ and $f(r_k\setminus\{p_k\})\subset r_k$.
\item For $z=x+iy$,
\[
|f(z)-(z-i)|\le\eta(y)\quad\text{if } y>0,\qquad
|f(z)-(z+i)|\le\eta(-y)\quad\text{if } y<0 .
\]
In particular, $f(z)=z-i+O(e^{-2\Im z})$ as $\Im z\to+\infty$ and
$f(z)=z+i+O(e^{2\Im z})$ as $\Im z\to-\infty$.
\item $f=A_1^{-1}\circ\Nsin\circ A_1$. Consequently, $f$ has no finite
asymptotic values and $\Sing(f)=\PS(f)=\Pcrit(f)=\{c_k:k\in\Z\}$. The immediate
basin $U_k$ of $c_k$ is simply connected, $U_k=U_0+k\pi$,
$\deg(f|_{U_k})=3$, $\ell_k\subset U_k$, and $r_k\subset \mathcal{J}(f)$ for every $k$.
Moreover,
\[
U_k\subset S_k=\{z\in\C:k\pi<\Re z<(k+1)\pi\},
\]
$\partial U_k$ contains exactly two poles of $f$, namely $p_k$ and $p_{k+1}$,
both accessible from $U_k$, and $U_k$ has exactly two accesses to $\infty$.
Finally, $\mathcal{J}(f)\cap\C=\bigcup_kJ_k$ with $J_k$ the component of $\mathcal{J}(f)\cap\C$
containing $r_k$, and every Fatou component of $f$ is eventually mapped into
some $U_k$; in particular, $f$ has no wandering domains.
\end{enumerate}
\end{lemma}
\begin{proof}
	(a) The function \(\cot\) is \(\pi\)-periodic and real on the real axis, so \(\cot\overline{z}=\overline{\cot z}\).\vspace{1.2mm}
	
	(b) We have \(f'(z)=1-\csc^2z=-\cot^2z\), which vanishes if and only if	\(\cos z=0\), \textit{i.e.,} \(z=c_k\) for some \(k\). Since \(\cot c_k=0\), each \(c_k\) is fixed.
	Using \(\cot(\pi/2+k\pi+w)=-\tan w\) and \(\tan w=w+w^3/3+O(w^5)\), we obtain	\[
	\displaystyle
	f(c_k+w)=c_k+w-\tan w
	=c_k-\frac{w^3}{3}+O(w^5),\;w\to 0.
	\]
	Since the first nonzero term of $f(c_k+w)-c_k$ is $-w^3/3$,
	we obtain $\deg_{c_k}f=3$.\vspace{1.2mm}
	
	(c) The poles of \(\cot z\) are the zeros \(k\pi\) of \(\sin z\); they are simple because \(\sin'(k\pi)=\cos(k\pi)=\pm1\ne0\). We obtain the two formulas from the identities
	\[
	\begin{aligned}
		\cot\!\left(\frac{\pi}{2}+k\pi+it\right)
		&=-\tan(it)=-i\tanh t, && t\in\mathbb{R},\\
		\cot(k\pi+it)
		&=\cot(it)=-i\coth t, && t\in\mathbb{R}\setminus\{0\}.
	\end{aligned}
	\]
	
	(d) This is Lemma~\ref{lem:estimates}(iii), since \(f(z)-(z\mp i)=\cot z\pm i\).\vspace{1.2mm}
	
	(e) The identity \(f=A_1^{-1}\circ\Nsin\circ A_1\) follows from
	Lemma~\ref{lem:affine} with \(g=\sin\), \(h=\cos\), \(A=A_1\) and \(c=1\).
	\[
	\begin{tikzcd}[column sep=large, row sep=large]
		z \arrow[r, "A_1"] \arrow[d, "f"']
		& A_1(z) \arrow[d, "\Nsin"] \\
		f(z) \arrow[r, "A_1"']
		& A_1(f(z))=\Nsin(A_1(z))
	\end{tikzcd}
	\] The map
	\(A_1^{-1}(z)=z-\pi/2\) sends \(k\pi\) to \(c_{k-1}\), the pole \(\pi/2+k\pi\) of
	\(\Nsin\) to \(p_k\), the line \(\{k\pi+it\}\) to \(\ell_{k-1}\) and the line
	\(\{\pi/2+k\pi+it\}\) to \(r_k\). Hence all statements except the inclusion
	\(U_k\subset S_k\) follow from Lemma~\ref{lem:sin} and
	Lemma~\ref{lem:conjugacy}; in particular, the immediate basin of \((k+1)\pi\) for
	\(\Nsin\), whose boundary contains the poles \(\pi/2+k\pi\) and \(\pi/2+(k+1)\pi\),
	is mapped onto \(U_k\), whose boundary therefore contains \(p_k\) and \(p_{k+1}\).
	Since there are no finite asymptotic values, the critical values \(c_k\) are the
	only finite singular values; as \(\{c_k\}\) is closed in \(\C\) and consists of
	fixed points, 
	\begin{align*}
		\Sing(f)=\PS(f)=\Pcrit(f)=\{c_k:k\in\Z\}.
	\end{align*} Finally, \(U_k\) is a
	connected subset of \(\mathcal{F}(f)\) containing \(c_k\in S_k\) and disjoint from the lines
	\(r_k,r_{k+1}\subset \mathcal{J}(f)\), which bound \(S_k\); hence \(U_k\subset S_k\).
\end{proof}
\begin{proposition}\label{prop:petals}
Let $f=\Ncos$, let $M\ge2$ and let $P_\pm=\{z\in\C:\pm\Im z>M\}$. Then the
following hold.
\begin{enumerate}
\item $f$ is univalent on the half-plane $H=\{z:\Im z>M-\frac12\}$, and in
particular on $P_+$.
\item $\overline{P_+}\subset f(P_+)\subset\{w:\Im w>M-1-\eta(M)\}
\subset\{w:\Im w>0\}$.
\item $P_+$ is a repelling petal of $f$ at infinity. In
Definition~\ref{def:petal}\textup{(c)} one may take $g\equiv1$, $t_0=1$,
$c_1=\frac12$, $c_2=2$ and $\delta=\pi/4$.
\end{enumerate}
The same statements hold for $P_-$, with $H$ replaced by $\{z:\Im
z<-M+\frac12\}$ and the inequalities in \textup{(b)} reversed. In particular,
$f(P_+)\cap P_-=f(P_-)\cap P_+=\emptyset$ and $f(P_+)\cap f(P_-)=\emptyset$.
\end{proposition}

\begin{proof}
Put $\varepsilon(z)=f(z)-(z-i)=\cot z+i$. By Lemma~\ref{lem:estimates}, we obtain
\begin{equation}\label{eq:eps}
\begin{aligned}
|\varepsilon(z)|&\le\eta(\Im z)\le\eta(M)<0.04 &&\text{for }\Im z\ge M,\\
|\varepsilon'(z)|&=|\csc^2z|\le\sinh^{-2}(\Im z) &&\text{for }\Im z>0.
\end{aligned}
\end{equation}

(a) The half-plane $H$ is convex. For $z_1,z_2\in H$, we have
\[
f(z_2)-f(z_1)-(z_2-z_1)=\cot z_2-\cot z_1=-\int_{[z_1,z_2]}\csc^2\zeta\,d\zeta,
\]
and $\Im\zeta>M-\frac12\ge\frac32$ on the segment $[z_1,z_2]$. By
Lemma~\ref{lem:estimates}, we have $|\csc^2\zeta|\le\sinh^{-2}(3/2)<\frac14$ there, so
\begin{align*}
	|f(z_2)-f(z_1)|\ge\frac34|z_2-z_1|.
\end{align*} Hence, $f$ is injective on $H$.\vspace{1.2mm}

(b) If $z\in P_+$, then 
\begin{align*}
	\Im f(z)=\Im z-1+\Im\varepsilon(z)>M-1-\eta(M),
\end{align*} since
$\eta$ is decreasing; and $M-1-\eta(M)\ge1-0.04>0$. To prove the first
inclusion, fix $w$ with $\Im w\ge M$, let $E=\{z:\Im z\ge M+\frac12\}$ and
define $\Phi(z)=w+i-\varepsilon(z)$. For $z\in E$ we have
\[
\Im\Phi(z)\ge M+1-\eta\left(M+\frac12\right)\ge M+\frac12,
\]
so $\Phi(E)\subset E$, and for $z_1,z_2\in E$,
\[
|\Phi(z_1)-\Phi(z_2)|=|\cot z_1-\cot z_2|\le\sinh^{-2}\left(M+\frac12\right)\,|z_1-z_2|
\le\frac{1}{30}|z_1-z_2|,
\]
by \eqref{eq:eps} and the convexity of $E$. Since $E$ is closed, the Banach fixed point theorem gives $z^*\in E$ with $z^*=w+i-\varepsilon(z^*)$, that is, $f(z^*)=z^*-i+\varepsilon(z^*)=w$. As $z^*\in E\subset P_+$, we obtain
$w\in f(P_+)$.\vspace{1.2mm}

(c) By (a) and the open mapping theorem, $f|_H$ is a homeomorphism onto the open set $f(H)$, so $G=(f|_H)^{-1}\colon f(H)\to H$ is holomorphic. Let $P'=f(P_+)$. Then $P'$ is a simply connected domain, being a conformal image of $P_+$, and it is unbounded because $P'\supset\overline{P_+}$ by (b).\vspace{1.2mm}

We claim that $\overline{P'}=f(\overline{P_+})$. The inclusion $\supset$ holds
by continuity. Conversely, if $w_n=f(z_n)\to w$ with $z_n\in P_+$, then
$|z_n-w_n-i|=|\varepsilon(z_n)|\le\eta(M)$, so $(z_n)$ is bounded and a
subsequence converges to some $z\in\overline{P_+}$ with $f(z)=w$. Since
$\overline{P_+}\subset H$, the map $G$ is defined and continuous on
$\overline{P'}\subset f(H)$, holomorphic on $P'$, and $G(P')=P_+$. For
$w\in\overline{P'}$ we have $G(w)\in\overline{P_+}$ and
$G(w)-w=G(w)-f(G(w))=i-\varepsilon(G(w))$; hence
\begin{equation}\label{eq:G}
|G(w)-w-i|\le\eta(\Im G(w))\le\eta(M)<0.04,\;
\Im G(w)>\Im w+\frac12 .
\end{equation}
We now verify conditions (a)--(c) of Definition~\ref{def:petal} for $G$ on $P'$.

\noindent {Condition (a).} We see that $\overline{G(P')}=\overline{P_+}\subset P'$ by (b).\vspace{1.2mm}

\noindent{Condition (b).} Every $w\in P'$ satisfies $\Im w>0$ by (b), so
\eqref{eq:G} gives $\Im G^n(w)>n/2$ for all $n\ge1$. Hence no point lies in $G^n(P')$ for all $n$, \textit{i.e.,} 
\begin{align*}
	\bigcap_{n\ge0}G^n(P')=\emptyset.
\end{align*}

\noindent{Condition (c).} Let $K\subset P'$ be compact, put
$R_K=\max_{z\in K}|\Re z|$, take $g\equiv1$ on $(t_0,+\infty)$ with $t_0=1$, $c_1=1/2$, $c_2=2$ and $\delta=\pi/4$, and let
\[
n_0=\Bigl\lceil\frac{2(R_K+1)}{\tan(\pi/8)}\Bigr\rceil .
\]
Let $z_0\in K$ and $z_n=G^n(z_0)$. For $n\ge1$ we have $z_n\in P_+$, so
$|z_n|>M>t_0$, and by \eqref{eq:G},
\begin{align*}
	\frac12<1-\eta(M)\le|G(z_n)-z_n|\le1+\eta(M)<2.
\end{align*} Moreover,
\begin{align*}
	\Im z_n>M+\frac{n-1}{2}\; \mbox{for}\; n\ge1
\end{align*} and, by \eqref{eq:G} and
Lemma~\ref{lem:estimates}(iv), we obtain
\[
|\Re z_n-\Re z_{n-1}|=|\Re\varepsilon(z_n)|\le\eta(\Im z_n)
\le4e^{-2M-(n-1)} .
\]
Summing over $n$ gives
\[
|\Re z_n-\Re z_0|\le\frac{4e^{-2M}}{1-e^{-1}}<7e^{-4}<1,
\]
so $|\Re z_n|<R_K+1$, while $\Im z_n>n/2$ because $\Im z_0>0$. Therefore, for $n\ge n_0$, the angle between $z_n$ and the positive imaginary axis is at most $\arctan\bigl(2(R_K+1)/n\bigr)\le\pi/8$. On the other hand, by \eqref{eq:G} the angle between $G(z_n)-z_n$ and $i$ is at most $\arcsin\eta(M)<\pi/8$. Consequently, the angle between $G(z)-z$ and $z$ is less than $\pi/4=\delta$ for every $z\in\bigcup_{n\ge n_0}G^n(K)$.\vspace{1.2mm}

Thus $P'$ is an attracting petal of $G$ at infinity, and since $f$ is
univalent on $P_+$, the half-plane $P_+$ is a repelling petal of $f$ at
infinity. Finally, complex conjugation $\sigma(z)=\overline{z}$ maps $P_+$ onto $P_-$, satisfies $f\circ\sigma=\sigma\circ f$ by Lemma~\ref{lem:cos}(a), and preserves moduli and angles between vectors; hence all statements carry over to $P_-$. The disjointness statements follow from (b) and its analogue for $P_-$, since
$f(P_+)\subset\{\Im w>0\}$ and $f(P_-)\subset\{\Im w<0\}$.
\end{proof}

\begin{remark}\label{rem:petals}
(i) Proposition~\ref{prop:petals} is the analogue for $\Ncos$ of
\cite[Example~4.5]{BFJK25}, where the half-planes $\{\pm\Im z>M\}$ are shown to be repelling petals of $\Nsin$ at infinity for sufficiently large $M$. We have verified Definition~\ref{def:petal} directly instead of using the sufficient condition of \cite[Proposition~4.4]{BFJK25}; this gives the explicit bound $M\ge2$. If one prefers to apply \cite[Proposition~4.4]{BFJK25}, then $d=1$, $a=-i$ for $P_+$ and $a=i$ for $P_-$, and, since repelling petals correspond to odd indices, $j=1$ in both cases.\vspace{1.2mm}

(ii) Since $\Nsin(z)-(z-i)=i-\tan z$ and $\Nsin'(z)-1=-\sec^2z$,
Lemma~\ref{lem:estimates} provides exactly the same estimates for $\Nsin$ as \eqref{eq:eps} does for $\Ncos$, and $\Nsin$ also commutes with complex conjugation. Hence the proof of Proposition~\ref{prop:petals} applies verbatim to $\Nsin$: for every $M\ge2$, the half-planes $\{\pm\Im z>M\}$ are repelling petals of $\Nsin$ at infinity.
\end{remark}

\section{Proofs of the main results}\label{sec:proofs}

\begin{proof}[\bf Proof of Theorem~\ref{thm:conj}]
We apply Lemma~\ref{lem:affine} with $g=\sin$. Since
\[
\sin(A_1(z))=\cos z,\; \sin(A_2(z))=\sin(iz)=i\sinh z,\;
\sin(A_3(z))=\cos(iz)=\cosh z,
\]
we obtain the first three identities in \eqref{eq:conj} (with $c=1$, $c=i$ and $c=1$, respectively). Since $\cos(A_2(z))=\cosh z$, Lemma~\ref{lem:affine} with $g=\cos$ gives $\Ncosh=A_2^{-1}\circ\Ncos\circ A_2$.\vspace{1.2mm}

For the reader's convenience, we also verify \eqref{eq:conj} directly. Using $\tan(z+\pi/2)=-\cot z$,
\[
A_1^{-1}\bigl(\Nsin(A_1(z))\bigr)
=\Bigl(z+\frac{\pi}{2}-\tan\Bigl(z+\frac{\pi}{2}\Bigr)\Bigr)-\frac{\pi}{2}
=z+\cot z=\Ncos(z).
\]
In view of $\tan(iz)=i\tanh z$ and $A_2^{-1}(w)=-iw$, we obtain
\[
A_2^{-1}\bigl(\Nsin(A_2(z))\bigr)=-i\bigl(iz-i\tanh z\bigr)=z-\tanh z
=\Nsinh(z).
\]
Using $\cot(iz)=-i\coth z$, we obtain
\[
A_2^{-1}\bigl(\Ncos(A_2(z))\bigr)=-i\bigl(iz-i\coth z\bigr)=z-\coth z
=\Ncosh(z),
\]
and $A_3=A_1\circ A_2$ gives $A_3^{-1}\circ\Nsin\circ A_3=A_2^{-1}\circ
\bigl(A_1^{-1}\circ\Nsin\circ A_1\bigr)\circ A_2=A_2^{-1}\circ\Ncos\circ
A_2=\Ncosh$. The remaining statements follow from
Lemma~\ref{lem:conjugacy}\textup{(b), (c)}.
\end{proof}
\begin{remark}
	We summarize the four conjugacies in the following diagram.
	\begin{equation}\label{eq:newton-conjugacy-diagram}
		\begin{tikzcd}[column sep=2.6em, row sep=1.3em]
			z \arrow[r, "A_1"]
			& z+\dfrac{\pi}{2} \arrow[r, "\Nsin"]
			& \Nsin(A_1(z)) \arrow[r, "A_1^{-1}"]
			& \Ncos(z) \\
			z \arrow[r, "A_2"]
			& iz \arrow[r, "\Nsin"]
			& \Nsin(A_2(z)) \arrow[r, "A_2^{-1}"]
			& \Nsinh(z) \\
			z \arrow[r, "A_3"]
			& iz+\dfrac{\pi}{2} \arrow[r, "\Nsin"]
			& \Nsin(A_3(z)) \arrow[r, "A_3^{-1}"]
			& \Ncosh(z) \\
			z \arrow[r, "A_2"]
			& iz \arrow[r, "\Ncos"]
			& \Ncos(A_2(z)) \arrow[r, "A_2^{-1}"]
			& \Ncosh(z)
		\end{tikzcd}
	\end{equation}
	Since $A_3=A_1\circ A_2$, we obtain the third row of
	\eqref{eq:newton-conjugacy-diagram} by substituting the first
	conjugacy into the fourth. In fact, we can express as \begin{equation}\label{eq:cosh-conjugacy-diagram}
		\begin{tikzcd}[column sep=large, row sep=large]
			z \arrow[r, "A_2"] \arrow[d, "\Ncosh"']
			& A_2(z) \arrow[r, "A_1"] \arrow[d, "\Ncos"']
			& A_3(z) \arrow[d, "\Nsin"] \\
			\Ncosh(z) \arrow[r, "A_2"']
			& \Ncos(A_2(z)) \arrow[r, "A_1"']
			& \Nsin(A_3(z))
		\end{tikzcd}
	\end{equation}
	Since $A_3=A_1\circ A_2$, we read the left square and the outer
	rectangle of \eqref{eq:cosh-conjugacy-diagram} to obtain
	\[
	\begin{aligned}
		\Ncosh
		&=A_2^{-1}\circ\Ncos\circ A_2\\
		&=A_2^{-1}\circ
		\bigl(A_1^{-1}\circ\Nsin\circ A_1\bigr)\circ A_2\\
		&=A_3^{-1}\circ\Nsin\circ A_3.
	\end{aligned}
	\]
\end{remark}
\begin{proof}[\bf Proof of Theorem~\ref{thm:general}]
A number $\varphi$ with $e^{2i\varphi}=-A/B$ exists because $-A/B\ne0$, and
$C\ne0$. Since $e^{iT(z)}=e^{i\varphi}e^{\lambda z}$ and
$e^{-iT(z)}=e^{-i\varphi}e^{-\lambda z}$, we have
\[
C\sin(T(z))=\frac{C}{2i}\bigl(e^{i\varphi}e^{\lambda z}
-e^{-i\varphi}e^{-\lambda z}\bigr)
=Ae^{\lambda z}-Ae^{-2i\varphi}e^{-\lambda z}=Ae^{\lambda z}+Be^{-\lambda z},
\]
because $-Ae^{-2i\varphi}=-A\cdot(-B/A)=B$. Thus $\sin\circ T=C^{-1}g$, and Lemma~\ref{lem:affine}, applied with $\sin$, $g$, $T$ and $C^{-1}$ in the roles of $g$, $h$, $A$ and $c$, gives $N_g=T^{-1}\circ\Nsin\circ T$. The zeros of $g$ are the points $z$ with $\sin(T(z))=0$, that is, $T(z)=k\pi$ for some $k\in\Z$. The remaining statements follow from Lemma~\ref{lem:conjugacy}(b) and Lemma~\ref{lem:sin}(a), since $T^{-1}$ maps unbounded sets onto unbounded sets.
\end{proof}
\begin{remark}
	Since $g=C(\sin\circ T)$, we apply Lemma~\ref{lem:affine} with
	$\sin$ as the original function, $T$ as the affine change of
	variable, and $C$ as the nonzero scalar factor. We represent the
	resulting conjugacy by the commutative diagram
	\begin{equation}\label{eq:general-conjugacy-diagram}
		\begin{tikzcd}[column sep=large, row sep=large]
			z \arrow[r, "T"] \arrow[d, "N_g"']
			& T(z) \arrow[d, "\Nsin"] \\
			N_g(z) \arrow[r, "T"']
			& \Nsin(T(z)).
		\end{tikzcd}
	\end{equation}
	Reading the two paths in \eqref{eq:general-conjugacy-diagram}, we obtain
	\[
	T\circ N_g=\Nsin\circ T\; \mbox{and}\;
	N_g=T^{-1}\circ\Nsin\circ T.
	\]
	In particular, $T$ sends the zero $T^{-1}(k\pi)$ of $g$ to the
	zero $k\pi$ of $\sin z$. By Lemma~\ref{lem:conjugacy}(b), $T^{-1}$
	maps the immediate basin of $k\pi$ for $\Nsin$ onto the immediate
	basin of $T^{-1}(k\pi)$ for $N_g$.
\end{remark}
\begin{proof}[\bf Proof of Theorem~\ref{thm:lc}]
By Theorem~\ref{thm:C}, $J(\Nsin)$ is locally connected. By
Theorem~\ref{thm:conj}, the Julia sets of $\Ncos$, $\Nsinh$ and $\Ncosh$ are the images of $J(\Nsin)$ under the homeomorphisms $A_1^{-1}$, $A_2^{-1}$ and $A_3^{-1}$ of $\Ch$, and by Theorem~\ref{thm:general}, we have $J(N_g)=T^{-1}(J(\Nsin))$ for every $g$ as in that theorem. Hence all these Julia sets are locally connected by Lemma~\ref{lem:conjugacy}(d). Each of them is compact, and it is
connected because the maps are Newton maps of entire functions \cite{BFJK14}. Hence Lemma~\ref{lem:torhorst} shows that the boundary of every Fatou component is locally connected.
\end{proof}

\begin{proof}[\bf Proof of Theorem~\ref{thm:direct}]
Part (i) is Proposition~\ref{prop:petals}.

For part (ii), fix $k\in\Z$. By Lemma~\ref{lem:cos}(e), $U_k$ is a simply connected invariant attracting basin, so it suffices to check hypotheses (a)--(c) of Theorem~\ref{thm:A}.\vspace{1.2mm}

\noindent{(a).} By Lemma~\ref{lem:cos}(e), $f$ has no finite asymptotic values, so $\Pasym(f)=\emptyset$. Moreover, $\Pcrit(f)=\{c_j:j\in\Z\}$ has no accumulation points in $\C$, so $\Acc(\Pcrit(f))=\emptyset$. Hence the set in (a) is empty.\vspace{1.2mm}

\noindent{(b).} For $j\ne k$, the point $c_j$ lies in the open set $U_j$, which is disjoint from $U_k$; hence $c_j\notin\overline{U_k}$ and $\Pcrit(f)\cap\overline{U_k}=\{c_k\}$. With the compact set $L=\{c_k\}\subset U_k$, the set $\Pcrit(f)\cap(\overline{U_k}\setminus L)$ is empty, and (b) holds (compare Remark~\ref{rem:L}).\vspace{1.2mm}

\noindent{(c).} By Lemma~\ref{lem:cos}(e), $\overline{U_k}\subset
\overline{S_k}=\{k\pi\le\Re z\le(k+1)\pi\}$. Let $D$ be the disc centred at $0$ of radius $R=(|k|+1)\pi+M+1$. If $z\in\overline{U_k}\setminus D$, then $|\Re z|\le(|k|+1)\pi$ and therefore $|\Im z|\ge|z|-|\Re z|\ge M+1>M$, that is,
$z\in P_+\cup P_-$. By part (i), $P_+$ and $P_-$ are repelling petals of $f$ at infinity, and
\[
U_k\cap f(P_+)\cap P_-\subset f(P_+)\cap P_-=\emptyset\; \mbox{and}\;
U_k\cap f(P_-)\cap P_+\subset f(P_-)\cap P_+=\emptyset .
\]
Hence $U_k$ is tame at infinity.\vspace{1.2mm}

Theorem~\ref{thm:A} now shows that $\partial U_k$ is locally connected, and Theorem~\ref{thm:B} applies to $U_k$.
\end{proof}

\begin{proof}[\bf Proof of Corollary~\ref{cor:structure}]
For $f=\Nsin$, statements (a)--(e) are contained in Lemma~\ref{lem:sin} and the formulas following it, except for the fact that $U_k$ lies in the strip $\{\pi/2+(k-1)\pi<\Re z<\pi/2+k\pi\}$; this holds because $U_k$ is connected, contains $k\pi$ and is disjoint from the lines $r_{k-1},r_k\subset J(\Nsin)$ (Lemma~\ref{lem:sin}(a), (b)). For (f), the petals are given by Remark~\ref{rem:petals}(ii), and the proof of hypothesis (c) in Theorem~\ref{thm:direct} applies verbatim with this strip in place of $S_k$. Statement (g) follows from Theorem~\ref{thm:lc}.\vspace{1.2mm}

For $f=\Ncos$, statements (a)--(e) are Lemma~\ref{lem:cos}, (f) is
Theorem~\ref{thm:direct}, and (g) is Theorem~\ref{thm:lc}.\vspace{1.2mm}

For $f=\Nsinh$ and $f=\Ncosh$, we use Theorem~\ref{thm:conj} and
Lemma~\ref{lem:conjugacy}. Since $A_2^{-1}(w)=-iw$ and
$A_3^{-1}(w)=-i(w-\pi/2)$, we have
\begin{align*}
A_2^{-1}(k\pi)&=-k\pi i,\; A_2^{-1}\Bigl(\frac{\pi}{2}+k\pi\Bigr)
=-\Bigl(k+\frac12\Bigr)\pi i,\\
A_3^{-1}(k\pi)&=\Bigl(-k+\frac12\Bigr)\pi i, 
A_3^{-1}\Bigl(\frac{\pi}{2}+k\pi\Bigr)=-k\pi i,
\end{align*}
and 
\begin{align*}
	A_2^{-1}(\{\Re z=s\})=\{\Im z=-s\},\;
	A_3^{-1}(\{\Re z=s\})=\{\Im z=\pi/2-s\}.
\end{align*} Hence, as sets, the fixed points, poles and lines of $\Nsin$ in Lemma~\ref{lem:sin} are mapped onto the corresponding objects for $\Nsinh$ and $\Ncosh$ listed in Table~\ref{tab:maps} (after relabelling the indices). Therefore (a), (c), (d) and (e) follow from the case $f=\Nsin$ treated above and Lemma~\ref{lem:conjugacy}, while (b) holds because $\tanh$ and $\coth$ are $\pi i$-periodic. Statement (g) follows from Theorem~\ref{thm:lc}.\vspace{1.2mm}

It remains to prove (f) for $\Nsinh$ and $\Ncosh$. Let $\lambda\in\C$ with $|\lambda|=1$, let $R(z)=\lambda z$, and let $P$ be a repelling petal at infinity of a map $F$, with $G=F^{-1}$ on $F(P)$ as in
Definition~\ref{def:petal}. Then $R^{-1}(P)$ is a repelling petal at infinity of $R^{-1}\circ F\circ R$: univalence, closures and intersections are preserved by $R$, and for $\widetilde G=R^{-1}\circ G\circ R$ we have 
\begin{align*}
	|\widetilde G(z)-z|=|G(R(z))-R(z)|
\end{align*} and $|z|=|R(z)|$, while the angle between $\widetilde G(z)-z$ and $z$ equals the angle between $G(R(z))-R(z)$ and
$R(z)$. Applying this with $R=A_2$ to the petals $\{\pm\Im z>M\}$ of $\Nsin$ (Remark~\ref{rem:petals}(ii)) and of $\Ncos$ (Proposition~\ref{prop:petals}), and using 
\begin{align*}
	A_2^{-1}(\{\pm\Im w>M\})=\{\pm\Re z>M\}
\end{align*} together with $\Nsinh=A_2^{-1}\circ\Nsin\circ A_2$ and $\Ncosh=A_2^{-1}\circ\Ncos\circ A_2$, we see that the half-planes $\{\pm\Re z>M\}$ are repelling petals at infinity of
$\Nsinh$ and $\Ncosh$ for every $M\ge2$. Finally, $A_2^{-1}$ maps the basins of $\Nsin$ and $\Ncos$ onto those of $\Nsinh$ and $\Ncosh$, and it maps the discs centred at $0$ onto themselves; hence the tameness of the basins of $\Nsin$ and $\Ncos$ established above is transferred to the basins of $\Nsinh$ and $\Ncosh$.
\end{proof}
\begin{remark}
	We put $\widetilde P=R^{-1}(P)$ and
	$\widetilde F=R^{-1}\circ F\circ R$. The forward maps and their
	inverse branches fit into the commutative diagram
	\begin{equation}\label{eq:petal-rotation-diagram}
		\begin{tikzcd}[column sep=4em, row sep=3em]
			\widetilde P
			\arrow[r, "R"]
			\arrow[d, bend right=18, "\widetilde F"']
			&
			P
			\arrow[d, bend left=18, "F"]
			\\
			\widetilde F(\widetilde P)
			\arrow[r, "R"']
			\arrow[u, bend right=18, "\widetilde G"']
			&
			F(P)
			\arrow[u, bend left=18, "G"]
		\end{tikzcd}
	\end{equation}
	Reading the downward and upward arrows in
	\eqref{eq:petal-rotation-diagram}, we obtain
	\begin{equation}\label{eq:petal-forward-inverse-compositions}
		R\circ\widetilde F=F\circ R\; \mbox{and}\;
		\widetilde G=R^{-1}\circ G\circ R.
	\end{equation}
	
	We now take $R=A_2$, where $A_2(z)=iz$. Since
	$\Im(A_2(z))=\Re z$, the map $A_2$ sends the half-planes
	$\{\pm\Re z>M\}$ onto $\{\pm\Im w>M\}$, with the same choice
	of sign. The forward compositions in the two cases are
	\begin{equation}\label{eq:hyperbolic-petal-compositions}
		\begin{tikzcd}[column sep=2.6em, row sep=1.4em]
			z \arrow[r, "A_2"] & iz \arrow[r, "\Nsin"]
			& \Nsin(iz) \arrow[r, "A_2^{-1}"] & \Nsinh(z) \\
			z \arrow[r, "A_2"] & iz \arrow[r, "\Ncos"]
			& \Ncos(iz) \arrow[r, "A_2^{-1}"] & \Ncosh(z).
		\end{tikzcd}
	\end{equation}
	Thus \eqref{eq:petal-forward-inverse-compositions} transfers
	the inverse-branch estimates to the corresponding half-planes
	for $\Nsinh$ and $\Ncosh$.
\end{remark}
\section*{Declarations}
\noindent{\bf Acknowledgments:}
The authors would like to thank the reviewers for their helpful suggestions and comments to improve the results of the paper.\vspace{1.2mm}

\noindent\textbf{Funding.} No funding was received for this work.\vspace{1.2mm}

\noindent\textbf{Data Availability Statement.} No datasets were generated or analysed during the current study. The figures were produced with Python scripts, which are available from the author on reasonable request.\vspace{1.2mm}

\noindent\textbf{Conflict of Interest.} The authors declare that there is no conflict of interest.\vspace{1.2mm}

\noindent\textbf{Clinical Trial Number.} Not applicable.\vspace{1.2mm}

\noindent\textbf{Author Contribution.} 
Both authors contributed equally to the conceptualization, formal analysis, investigation, and writing of this manuscript. Both authors have read and agreed to the published version of the manuscript.

\end{document}